\documentclass[12pt]{article}
\usepackage{amsfonts}
\usepackage{amssymb}
\usepackage{amsmath}
\usepackage{color}
=0pt
\def\sqr#1#2{{\vcenter{\vbox{\hrule height.#2pt
              \hbox{\vrule width.#2pt height#1pt \kern#1pt \vrule
width.#2pt}
              \hrule height.#2pt}}}}
\def\signed #1{{\unskip\nobreak\hfil\penalty50
              \hskip2em\hbox{}\nobreak\hfil#1
              \parfillskip=0pt \finalhyphendemerits=0 \par}}
\def\endpf{\signed {$\sqr69$}}
\def\dbC{{\mathbb{C}}}

\def\dbF{{\mathbb{F}}}

\def\dbH{{\mathbb{H}}}

\def\dbN{{\mathbb{N}}}

\def\dbP{{\mathbb{P}}}

\def\dbR{{\mathbb{R}}}
\def\dbS{{\mathbb{S}}}

\def\a{\alpha}
\def\b{\beta}

\def\d{\delta}
\def\e{\varepsilon}

\def\l{\lambda}

\def\n{\nu}

\def\f{\varphi}
\def\th{\theta}
\def\o{\omega}

\def\3n{\negthinspace \negthinspace \negthinspace }
\def\2n{\negthinspace \negthinspace }
\def\1n{\negthinspace }
\def\ns{\noalign{\smallskip} }

\def\ds{\displaystyle}
\def\G{\Gamma}
\def\D{\Delta}

\def\L{\Lambda}
\def\Si{\Sigma}

\def\O{\Omega}
\def\cA{{\cal A}}

\def\cE{{\cal E}}
\def\cF{{\cal F}}

\def\cH{{\cal H}}

\def\cO{{\cal O}}

\def\cl{{\cal l}}
\def\mE{{\mathbb{E}}}

\def\no{\noindent}

\def\ss{\smallskip}
\def\ms{\medskip}
\def\bs{\bigskip}
\def\q{\quad}
\def\qq{\qquad}

\def\pa{\partial}

\def\wt{\widetilde}
\def\cd{\cdot}
\def\cds{\cdots}

\def\dist{\hbox{\rm dist$\,$}}

\def\supp{\hbox{\rm supp$\,$}}

\def\cl{\overline}

\def\({\Big (}
\def\|{\Big |}
\def\){\Big )}
\def\[{\Big[}
\def\]{\Big]}

\def\={\buildrel \triangle \over =}

\def\be{\begin{equation}}
\def\bel{\begin{equation}\label}
\def\ee{\end{equation}}
\def\bea{\begin{eqnarray}}
\def\eea{\end{eqnarray}}
\def\bt{\begin{theorem}}
\def\et{\end{theorem}}
\def\bc{\begin{corollary}}
\def\ec{\end{corollary}}
\def\bl{\begin{lemma}}
\def\el{\end{lemma}}
\def\bp{\begin{proposition}}
\def\ep{\end{proposition}}
\def\br{\begin{remark}}
\def\er{\end{remark}}
\def\ba{\begin{array}}
\def\ea{\end{array}}
\def\bd{\begin{definition}}
\def\ed{\end{definition}}

\newtheorem{lemma}{Lemma}[section]
\newtheorem{remark}{Remark}[section]

\newtheorem{theorem}{Theorem}[section]
\newtheorem{corollary}{Corollary}[section]

\newtheorem{definition}{Definition}[section]
\newtheorem{proposition}{Proposition}[section]
\newtheorem{condition}{Condition}[section]

\makeatletter
   
   \@addtoreset{equation}{section}
\makeatother

\begin{document}
\title{\bf  An Internal Observability  Estimate for Stochastic Hyperbolic Equations\thanks{This work is partially supported by the NSF of China
under grants 11231007, 11322110, 11371084 and
11471231, by
the Fundamental Research Funds for the Central
Universities under grant  2412015BJ011 and
2015SCU04A02, and by  the
Chang Jiang Scholars Program from Chinese
Education Ministry. The authors highly appreciate the anonymous referees for their constructive
comments which led to this improved version.}}

\author{Xiaoyu Fu\thanks{School of Mathematics, $\q$
Sichuan University, $\q$ Chengdu 610064, $\q$ China. $\q$
E-mail address:
xiaoyufu@scu.edu.cn.}\and Xu Liu\thanks{School of
Mathematics and Statistics, Northeast Normal
University, Changchun 130024,  China. E-mail
address: liux216@nenu.edu.cn.}\and Qi
L\"u\thanks{School of Mathematics,  Sichuan
University, Chengdu 610064,  China. E-mail
address: lu@scu.edu.cn.}\and Xu
Zhang\thanks{School of Mathematics, $\q$
Sichuan University, $\q$ Chengdu 610064, $\q$ China. $\q$
E-mail address: zhang\_xu@scu.edu.cn.}}

\date{}

\maketitle

\begin{abstract}
\no  This paper is addressed to establishing  an
internal observability estimate for some linear
stochastic hyperbolic equations.  The key is to
establish a new global Carleman estimate for
forward stochastic hyperbolic equations  in  the
$L^2$-space. Different from the deterministic
case, a delicate analysis on the adaptedness for
some stochastic processes is required in the
stochastic setting.
\end{abstract}

\bs

\no{\bf 2010 Mathematics Subject
Classification}. Primary 93B05; Secondary 93B07,
93C20.

\bs

\no{\bf Key Words}. Stochastic hyperbolic
equation, observability estimate, global Carleman
estimate, adaptedness, optimal control.

\section{Introduction and main result}

Let $T>0$ and  $(\Omega, \mathcal{F},
\{\mathcal{F}_t\}_{t\geq 0}, \mathcal{P})$ be a
complete filtered probability space,  on which a
one-dimensional standard Brownian motion
$\{B(t)\}_{t\geq 0}$ is defined such that
$\mathbb{F}=\{\mathcal{F}_t\}_{t\geq 0}$ is the
natural filtration generated by $B(\cdot)$,
augmented by all the $\mathcal{P}$-null sets in
$\mathcal{F}$. Let $\mathcal{H}$ be a Banach
space, and let $C([0,T];\mathcal{H})$ be the
Banach space of all $\mathcal{H}$-valued
strongly continuous functions defined
on $[0,T]$. We denote by
$L^2_{\mathbb{F}}(0,T;\mathcal{H})$ the Banach
space consisting of all $\mathcal{H}$-valued
$\dbF$-adapted processes $X(\cdot)$ such that
$\mathbb{E}(|X(\cdot)|^2_{L^2(0,T;\mathcal{H})})<\infty$,
with the canonical norm; by
$L^\infty_{\mathbb{F}}(0,T;\mathcal{H})$ the
Banach space consisting of all
$\mathcal{H}$-valued $\dbF$-adapted essentially
bounded processes; and by
$L^2_{\mathbb{F}}(\Omega; C([0,T];\mathcal{H}))$
the Banach space consisting of all
$\mathcal{H}$-valued $\dbF$-adapted continuous
processes $X(\cdot)$ such that
$\mathbb{E}(|X(\cdot)|^2_{C([0,T];\mathcal{H})})<\infty$.
 Similarly, one can define $L^\infty_\mathbb{F}(\Omega; C^m([0, T];
\mathcal{H}))$ for any positive integer $m$.

Let  $G \subset \mathbb{R}^{n}$ (for some $n\in\mathbb{N}$) be a  nonempty
bounded domain with a $C^{2}$ boundary $\G$.
Set $ Q=(0,T) \times G$ and  $\Si=(0,T) \times
\G$.  Assume that $ b^{ij} \in C^2(\cl{G})$
$(i,j = 1,2,\cdots, n)$
 satisfy
 \bel{gz1}
b^{ij}(x)= b^{ji}(x),\quad\forall\  x\in
\overline{G}, \ee
and for some constant $s_0
> 0$,
\begin{equation}\label{gz2}
 \sum_{i,j=1}^nb^{ij}(x)\xi^{i}\xi^{j} \geq s_0 |\xi|^2,
\,\,\,\,\,\,\,\,\,\, \forall \  (x,\xi)=
(x,\xi^{1}, \cdots, \xi^{n}) \in \overline{G}
\times \mathbb{R}^{n}.
\end{equation}

Let us consider the following stochastic hyperbolic
equation:
\begin{eqnarray}{\label{system1}}
\left\{
\begin{array}{lll}\ds
\ds dy_{t} -
\sum_{i,j=1}^{n}(b^{ij}y_{x_i})_{x_j}dt = \big(
b_1 y + f \big)dt + \big(b_2 y + g\big)dB(t)
&{\mbox {in }} Q,
 \\
\ns\ds  y = 0  &\mbox{on } \Si, \\[3mm]
\ns\ds  y(0) = y_0,\  y_{t}(0) = y_1 &\mbox{in }
G,
\end{array}
\right.
\end{eqnarray}
where  $(y_0, y_1) \in L^2(G)\times H^{-1}(G)$,
$f, g \in L^2_{\dbF}(0,T;L^2(G))$,
 \begin{equation} \label{aibi}
  b_1  \in L_{\dbF}^{\infty}(0,T;L^{p}(G))
\mbox{ with }  p\in [n,\infty] \q \mbox{ and }\q
 b_2 \in
L_{\dbF}^{\infty}(0,T;L^{\infty}(G)).
\end{equation}
Also,  set
$$ H_{T}=
L_{\dbF}^2 (\O; C([0,T];L^2(G)))\bigcap L_{\dbF}^2
(\O; C^{1}([0,T];H^{-1}(G)))
$$
and
$$
\cH_T=L^2_\dbF(\O;C^1([0,T];L^2(G)))\bigcap
L^2_\dbF(\O;C([0,T];H_0^1(G))).
$$
Then $H_T$ and $\cH_T$ are Banach spaces with
the canonical norms. In this paper, we use the following notion of solution for the equation \eqref{system1}.
\begin{definition} \label{def solution to sys}
A function $y\in H_T$ is called  a solution to
the equation \eqref{system1},  if the following
conditions
hold: \\[3mm]$(1)$ $y(0) = y_0$  and
$y_t(0) = y_1$ in $G$, $\mathcal{P}$-a.s.  \\[3mm]
$(2)$ For any $t \in (0,T)$ and $\f \in
C^2(\overline{G})\cap C_0^1(G)$,  it holds that
\begin{equation} \label{solution to sysh}
\begin{array}{ll}\ds
  \langle
y_t(t),\f\rangle_{H^{-1}(G),H_0^1(G)} - \langle
y_t(0),\f\rangle_{H^{-1}(G),H_0^1(G)}
\\ \ns\ds= \int_0^t \int_G \Big[
\sum_{i,j=1}^n
\big(b^{ij}(x)\f_{x_i}(x)\big)_{x_j} y(s, x) +
\Big(b_1(s,x) y(s,x)  + f(s,x)
\Big)\f(x)\Big]dxds \\
\ns\ds \q   + \int_0^t \int_G \Big(b_2(s,x)
y(s,x) + g(s,x)\Big)\f(x) dxdB(s),\q
\mathcal{P}\mbox{-a.s.}
\end{array}
\end{equation}
\end{definition}

For any initial value
$(y_0, y_1) \in L^2(G) \times H^{-1}(G)$, it is easy to show that the
equation (\ref{system1}) admits a unique
solution $y \in H_T $.

Let  $\G_0$ be a part of the boundary of $G$ satisfying
certain conditions, which
will be specified  later. For any given
constant $\d>0$, set
$$
 \cO_\d(\G_0) = \Big\{ x\in G ;\,
\dist(x,\G_0)<\d \Big\}.
$$
Write
\begin{equation}\label{r1r2}
r_1 = |b_1|_{L_{\dbF}^{\infty}(0,T;L^{p}(G))} \q
\mbox{ and }\q r_2 =
|b_2|_{L_{\dbF}^{\infty}(0,T;L^{\infty}(G))}.
\end{equation}
The main purpose of this paper is to establish
the following inequality:
\begin{equation} \label{inobser esti2.1}
\begin{array}{ll}\ds
|(y_0,y_1)|_{L^2(G)\times H^{-1}(G)}
\\ \ns\ds \leq Ce^{C\big(r_1^{\frac{1}{2-n/p}}+r_2^2\big)}
\Big(|y|_{L^2_{\mathbb{F}}(0, T;
L^2(\mathcal{O}_\delta(\Gamma_0)))}+|(f,
g)|_{(L^2_\mathbb{F}(0, T; L^2(G)))^2}\Big),
\end{array}
\end{equation}
where $y$ is the solution to (\ref{system1})
corresponding  to any given initial value $(y_0,
y_1)$.  Here and henceforth,
$C$ denotes a generic positive constant (which may be different from line to line), depending only
on $G$, $T$, $\G_0$, $b^{ij}$ ($i,j=1,\cds,n$),
$\d$,   and $d(\cdot)$ and  $\mu_0$ in
Condition \ref{condition of d} (to be given
later).

The inequality  \eqref{inobser esti2.1}  is
called an observability estimate for
(\ref{system1}).  For the case that $(f,g)=0$ in (\ref{system1}), this inequality means that the initial  energy of
a solution  in  the time $t=0$ can be bounded
by  its partial energy in  the local observation domain
$\cO_{\d}(\Gamma_0)$ in the time duration $[0,T]$.
Such kind of inequalities are closely  related
to control and state observation problems of
deterministic/stochastic hyperbolic equations.
For example, they can be applied to  a study  of
the controllability (e.g.
\cite{Bardos-Lebeau-Rauch1, Coron, Fu-Yong-Zhang1,
Lions1, Zh2, Zua}) and also inverse problems (e.g.
\cite{Kli2,Lav})
 for deterministic hyperbolic equations.  There exist numerous works devoted
to observability estimates for deterministic
hyperbolic equations. However,  there are only a very few  works addressed to similar problems but for stochastic
hyperbolic equations (\cite{Lu1, LZ, Zh3}).

Up to now,  there are several methods to
establish observability estimates for
deterministic hyperbolic equations, such as the (Rellich-type)
multiplier method (\cite{Lions1}), the
non-harmonic Fourier series technique
(\cite{Russell}),  the method of micro-local
analysis (\cite{Bardos-Lebeau-Rauch1}) and the method of
global Carleman estimate (\cite{Zh2}). The
multiplier  method is only applicable to some
very special hyperbolic equations.  Indeed, even
for a deterministic  hyperbolic equation,
 the observability estimate cannot be derived by this method directly for the case that the coefficients of lower
 order terms depend on both the space variable and time variable.
Also,  the non-harmonic Fourier series technique
has  restrictions not only on the coefficients,
but also on the  spatial  domain $G$ (requiring the domain $G$ to have some special shapes).
Furthermore,  since  the propagation of
singularities for stochastic partial
differential  equations is far from being well-understood, how to use the method of micro-local
analysis in the stochastic framework to
establish observability estimates is still unclear. Therefore, the Carleman estimate method
turns out to be a useful tool to establish
 observability estimates for stochastic
hyperbolic equations.

In \cite{Zh3}, by means of
a global Carleman estimate, a boundary
observability estimate for the equation
\eqref{system1} (with $(b^{ij})_{1\leq i,j \leq
n}=I_n$, the $n\times n$ identity matrix) was obtained:
\begin{equation} \label{obser esti22}
\begin{array}{ll}\ds
|(y(T),y_t(T))|_{L^2_{\cF_T}(\O;
H_0^1(G)\times L^2(G))}
\\ \ns\ds \leq Ce^{C(r_1^2 +
r^2_2)} \Big(\Big|\frac{\partial
y}{\partial \nu}\Big
|_{L^2_{\dbF}(0,T;L^2(\G_0))} + |(f,
g)|_{(L^2_{\dbF}(0,T;L^2(G)))^2}\Big),
\end{array}
\end{equation}
where $y$ solved the equation \eqref{system1}
associated to an initial data $(y_0,y_1)\in
H_0^1(G)\times L^2(G)$, and $\nu=\nu(x)=(\nu^1,\nu^2,\cdots,\nu^n)$ denotes
the unit outward normal vector of $\Omega$ at $x\in \Gamma$.  Also, in (\ref{obser esti22}), $T$ was
required to satisfy the condition:
$$
\frac{4+5c}{9c} \min_{x\in\overline
G}|x-x_0|^2> c^2 T^2 >
4\max_{x\in\overline G}|x-x_0|^2,
$$
for some $c\in (0, 1)$ and $x_0\in
\mathbb{R}^n\setminus\overline G$.

In \cite{Lu1},  by virtue of another global Carleman
estimate,  the result in \cite{Zh3} was improved
to the following boundary observability inequality:
\begin{eqnarray}\label{1**}
\begin{array}{ll}\ds
 |(y_0,y_1)|_{H_0^1(G)\times L^2(G)}
\\ \ns\ds  \leq C e^{C\big(r_1^{\frac{1}{3/2-n/p}}+ r_2^2\big)}\Big(\Big|\frac{\partial y}{\partial \nu}\Big |_{L^2_{\dbF}(0,T;L^2(\G_0))} + |(f, g)|_{(L^2_{\dbF}(0,T;L^2(G)))^2}\Big),
\end{array}
\end{eqnarray}
with $T >
2\max\limits_{x\in\overline G}|x-x_0|$  \Big(for the case that $(b^{ij})_{1\leq i,j \leq
n}=I_n$\Big).  Notice that in $(\ref{1**})$, the power  of $r_1$  is smaller than that in  (\ref{obser esti22}) (Indeed,
$\frac{1}{3/2-n/p}\leq 2$).
Also, an internal observability
estimate was established in \cite{Lu1}:
\begin{equation} \label{inobser esti2.3}
\begin{array}{ll}\ds
|(y_0,y_1)|_{H_0^1(G)\times L^2(G)}
\\
\ns\ds \leq C e^{C\big(r_1^{\frac{1}{3/2-n/p}}+
r_2^2\big)}\Big(\big|\nabla y\big
|_{L^2_{\dbF}(0,T;L^2(\cO_\d(\G_0)))} + |(f,
g)|_{(L^2_{\dbF}(0,T;L^2(G)))^2}\Big).
\end{array}
\end{equation}

The main difference between \eqref{inobser
esti2.1} and \eqref{inobser esti2.3} is that the
inequality \eqref{inobser esti2.3} provides an
observability estimate of the $H^1$-norm for
solutions to the equation \eqref{system1},  but
the inequality \eqref{inobser esti2.1} is an
estimate of the $L^2$-norm type. Compared with
the known inequality \eqref{inobser esti2.3}, the estimate
\eqref{inobser esti2.1} has more applications. For example, one application of \eqref{inobser esti2.1} is the
stabilization of stochastic hyperbolic
equations (but the detailed analysis is beyond the scope of this paper
and will be presented in our forthcoming work).
On the other hand, the inequality \eqref{inobser
esti2.1}  can be used to solve state
observation problem (See \cite{Lu1} for
example).

It is not clear whether it is possible to derive \eqref{inobser
esti2.1} from \eqref{inobser esti2.3} directly.
Indeed, this turns out to be  very difficult,  even for deterministic linear hyperbolic equations with some lower order terms  (See
\cite{Fu-Yong-Zhang1, I}). Also, compared with the
deterministic problem, the stochastic setting
will bring this problem some new difficulties. Actually, as we shall see later, a
suitable auxiliary optimal control problem (different from the deterministic context)
has to be constructed to guarantee the adaptedness of the related
stochastic processes.

Before giving our main result,   let us first
introduce some assumptions on $(b^{ij})_{1\leq
i,j\leq n}$ $(i, j=1, \cdots, n)$ and $T$.

\begin{condition}
\label{condition of d} There exists a positive
function $d(\cdot) \in
C^2(\overline{G})$ with the property
that \linebreak $\min\limits_{x\in \overline{G} }|\nabla d(x)| >
0$ such that, for some constant $\mu_0>0$,  the interior compatibility condition and the boundary compatibility condition are
satisfied:
$$\ds
\sum_{i,j=1}^n \sum_{i',j'=1}^n\!\!\Big(
2b^{ij'}(b^{i'j}d_{x_{i'}})_{x_{j'}} -
b^{ij}_{x_{j'}}b^{i'j'}d_{x_{i'}} \!\Big)\xi^{i}\xi^{j} \geq \mu_0 \sum_{i,j=1}^n
b^{ij}\xi^{i}\xi^{j}, \quad \forall\
(x,\xi^{1},\cdots,\xi^{n}) \in \overline{G}
\times \mathbb{R}^n.
$$
\end{condition}

In the sequel, we shall choose the set $\G_0$ as follows:
\begin{eqnarray}\label{def gamma0}
\G_0= \Big\{ x\in \G;
\sum_{i,j=1}^nb^{ij}(x)d_{x_i}(x)\nu^{j}(x) >0
\Big\},
\end{eqnarray}
where the function $d(\cd)$ is given in
Condition $\ref{condition of d}$. Also, write
\bel{20160115e2}
G_0=\cO_\d(\G_0)\quad \mbox{ and
}\quad\Sigma_0=\Gamma_0\times(0, T).
\ee

\begin{remark}
Notice that Condition $1.1$ is a sufficient condition for establishing Carleman estimates for deterministic linear hyperbolic operators $\ds \pa_t^2-\sum_{i,j=1}^n\pa_{x_j}(b^{ij}\pa_{x_i})$.
If $(b^{ij})_{1\leq i,j\leq n}=I_n$, then $d(x)=|x-x_0|^2$ satisfies
Condition $1.1$ with $x_0$ being
any given point in
$\dbR^n\setminus\overline G$. On the other hand, Condition $1.1$ can  also  be regarded as a special case of the pseudo-convexity condition in \cite{H}.  In fact, for the wave operator $\pa_t^2-\D$, if we set $a(x,\xi)=|\xi|^2$ and $d(x)=|x-x_0|^2$, then it is easy to check that
 $$
 \{a,\{a, d\}\}(x,\xi)=4|\xi|^2>0,\q  \forall\
(x,\xi) \in \overline{G}
\times (\dbR^n\setminus \{0\}),
 $$
 where $\{a,d\}$ denotes the Poission
bracket of $a$ and $d$, i.e.,
  $$\ba{ll}\ds
  \{a,d\}(x,\xi)=\sum_{j=1}^n\(\frac{\pa a}{\pa\xi_j}\cdot \frac{\pa d}{\pa x_j}-\frac{\pa a}{\pa x_j}\cdot\frac{\pa d}{\pa\xi_j}\),\q \forall\; (x,\xi)\in\cl{G}\times\dbR^n\setminus\{0\}.
  \ea$$
Moreover, it is easy to see that there is no critical point of the function $d(\cdot)$ in $\overline G$.
\end{remark}

\begin{remark}\label{rm1}
In \cite{Fu-Yong-Zhang1},  Condition $\ref{condition of d}$ was used
to establish  an internal observability estimate
for deterministic  hyperbolic equations. We refer to
\cite{Fu-Yong-Zhang1, liu} for more explanation on
Condition $\ref{condition of d}$ and some
interesting nontrivial examples. Also, as \cite[Example 1.1]{liu} shows, there exists some example for which Condition $\ref{condition of d}$ fails.
\end{remark}

In what follows,  set
 \bel{20160115e1}
M_0= \min\limits_{x\in\overline{G}} \sum_{i,j=1}^nb^{ij}d_{x_i}d_{x_j},\q  M_1= \max\limits_{x\in\overline{G}} \sum_{i,j=1}^nb^{ij}d_{x_i}d_{x_j} \q \mbox{ and }\q d_0=\max\limits_{x\in\G}\sum_{i,j=1}^n b^{ij}d_{x_i}\nu^j,
 \ee
and   define\footnote{From the proof of \cite[Theorem 5.1]{Fu-Yong-Zhang1}, it is easy to see that the number $d_0$ defined in \eqref{20160115e1} is positive. Hence, the set $\G_0$ given by \eqref{def gamma0} is not empty.}
\begin{equation}\label{fllz101}
T_0=\max\Big\{2\sqrt{M_1},  1+\frac{24\sqrt{n}d_0}{\min\{1,s_0\}}\(1+\frac{1}{s_0^{3/2}} \sum\limits_{i, j=1}^{n}|b^{i j}|_{C(\overline{G})}+\frac{1}{s_0}\)\Big\},
\end{equation}
where $s_0$ is the constant appeared in $(\ref{gz2})$.

\begin{remark}
It is easy to check that if $d(\cdot) $
satisfies Condition $\ref{condition of d}$,
then for any   $a \geq 1$ and $b \in
\mathbb{R}$, the function $\tilde{d}(\cdot) =
ad(\cdot)+b$ still satisfies this  condition
when  $\mu_0$ is replaced by $a\mu_0$.  Therefore,  throughout this paper, we may assume that $d(\cdot)$ and $\mu_0$ satisfy that
 \bel{t0}\ds
 \mu_0>\ds\frac{9T^2_0}{M_0}\q \mbox{ and }\q
M_0\geq \max\limits_{x\in\overline{G}} d(x).
 \ee
 \end{remark}

The main result of this paper is stated
as follows.
\begin{theorem}\label{inobser}
Assume that Condition $\ref{condition of
d}$  holds. Then, for any $T>T_0$ (defined by \eqref{fllz101}), the
observability inequality \eqref{inobser esti2.1}
holds for any solution to the equation
$(\ref{system1})$.
\end{theorem}
\br
The restriction on T in Theorem $\ref{inobser}$ is a technical condition, and $T_0$ is not sharp.  However, this condition plays a key role in our proof of Theorem $\ref{inobser}$.  It is reasonable to
expect that it can be improved to a better
one as that in \cite{Lu1} (for the estimates \eqref{1**} and \eqref{inobser esti2.3}), but this is an unsolved problem.
 \er
 \br
The condition $(\ref{t0})$ is  relevant to the interior behavior/property of the diffusion,  and it will  play a key role in the estimates on the energy terms $($see  $(2.13)$-$(2.15)$  in the proof of Theorem $\ref{5.24-lm1})$.  On the other hand, the  assumption  on the time $T$ in Theorem $\ref{inobser}$  is  relevant to the diffusion/reflection on the boundary. This assumption will play a key role in the estimates on the  boundary term $($see  Step $4$ in the proof of Theorem $\ref{5.24-lm1})$.   If  one considers  a special case, i.e. $(b^{ij})_{1\leq i, j\leq n}=I_n$, then $s_0=1$ and  we take $d(x)=|x-x_0|^2$, the corresponding condition on $T$ is  the following:
 $$
T>T_0=\max\Big\{4\max\limits_{x\in \overline{G}}|x-x_0|, \;\;1+{48\sqrt{n}}(n+2)\max\limits_{x\in\Gamma}\[(x-x_0)\cdot\nu(x)\]\Big\}.
 $$
 \er

The rest of this paper is organized as follows. In Section \ref{sss2}, we present a key weighted identity for partial differential operators of second order with symmetric coefficients.
Section \ref{sec-H1} is devoted to establishing
a Carleman estimate for  deterministic
hyperbolic equations in the $H^1$-space. In
Section \ref{ss6}, an auxiliary optimal control problem is introduced and analyzed. In
Section \ref{ss7}, a global
Carleman estimate for stochastic hyperbolic
equations  in the $L^2$-space is derived. In
Section \ref{sec-en-back},  energy estimates for random
hyperbolic equations and backward stochastic
hyperbolic equations are given. Section \ref{sec-pr in} is
devoted to a proof of our main result (i.e., Theorem
\ref{inobser}).  Finally, in Appendices A and B, we give
the  proofs of  some technical results.

\section{A weighted identity for  partial differential operators of second order with symmetric coefficients}\label{sss2}

In this section, we show a  pointwise
weighted identity for partial differential operators of second order
with symmetric coefficients, which will play a
crucial role in the sequel.

\bl\label{c1t1}
Assume that
 $a^{ij}=a^{ji}\in C^2(\dbR^m)$  ($i,j=1,2,\cdots,m$) for some $m\in\dbN$, $u,\ \ell\in C^2(\dbR^m)$ and  $\Psi\in C^1(\dbR^m)$.
 Set $\th=e^{\ell }$ and $v=\th u$. Then
 \bel{c1e2a}
 \3n\ba{ll}
 &\displaystyle
 \th^2\left|\sum_{i,j=1}^m(a^{ij}u_{x_i})_{x_j}\right|^2+ 2\sum_{j=1}^m\left[\sum_{i,i',j'=1}^m\(2a^{ij} a^{i'j'}\ell_{x_{i'}}v_{x_i}v_{x_{j'}}
 -a^{ij}a^{i'j'}\ell_{x_i}v_{x_{i'}}v_{x_{j'}}\)\right.\\
 \noalign{\ss}
 &\displaystyle\left.\qq\qq\qq\qq\qq\qq\qq\qq\qq\qq+\Psi \sum_{i=1}^m
 a^{ij}v_{x_i}v-\Lambda\sum_{i=1}^m a^{ij} \ell_{x_i}v^2\right]_{x_j}\\
 \noalign{\ss}
 &\displaystyle
\ge 2\sum_{i,j=1}^m
\Big[\sum_{i',j'=1}^m\(2a^{ij'}\big(a^{i'j}\ell_{x_{i'}}\big)_{x_{j'}}
-
 \big(a^{ij}a^{i'j'}\ell_{x_{i'}}\big)_{x_{j'}}\)+\Psi a^{ij}
 \Big]v_{x_i}v_{x_j}\\
 \ns&\ds\q
 +2\sum_{i,j=1}^m a^{ij}\Psi_{x_j}vv_{x_i}+Bv^2,
 \ea
 \ee
where
\bel{c1e3a}
 \Lambda=-\sum_{i,j=1}^m
 (a^{ij}\ell_{x_i}\ell_{x_j}-a^{ij}_{x_j}\ell_{x_i}
 -a^{ij}\ell_{x_ix_j})-\Psi,\qq   B=2\Lambda\Psi-2
 \sum_{i,j=1}^m\big(\Lambda a^{ij}\ell_{x_i}\big)_{x_j}.
\ee
\el
 \br
Lemma $\ref{c1t1}$  looks very similar to
Theorem $4.1$ in $\cite{Fu-Yong-Zhang1}$. The only
difference is about the regularity on the auxiliary
function $\Psi$. In $\cite{Fu-Yong-Zhang1}$, $\Psi$ was required to be
in $C^2(\dbR^m)$. But here we weaken  this
requirement to be $\Psi\in C^1(\dbR^m)$.
Note that,  the  choice of $\Psi$ usually depends
on  coefficients of the principal operator under consideration $($See
the equation  $(4.13)$ in $\cite{Fu-Yong-Zhang1})$.
Hence,  this implies that we only need the
$C^2$-regularity for coefficients of principal
operators, rather than the $C^3$-regularity required in
 $\cite{Fu-Yong-Zhang1}$.
 \er

\noindent{\bf Proof of Lemma \ref{c1t1}.}
Recalling $\th=e^\ell$ and $v=\th u$, we see
that
  $\th u_{x_i}=v_{x_i}-\ell_{x_i}v\ (i=1,2,\cds,m). $ Proceeding exactly as \cite[Theorem 4.1]{Fu-Yong-Zhang1},  we obtain that
 $$
 \displaystyle
-\th
\sum_{i,j=1}^m(a^{ij}u_{x_i})_{x_j}=I_1+I_2,
 $$
where
\bel{2c2t4}
 I_1=-\sum_{i,j=1}^m (a^{ij}v_{x_i})_{x_j}+\L v,\qq
 I_2=2\sum_{i,j=1}^m a^{ij}\ell_{x_i}v_{x_j}+\Psi v.
 \ee
This implies that
 \bel{c1e5}
 \ba{ll}
 \displaystyle\th^2\Big|\sum_{i,j=1}^m(a^{ij}u_{x_i})_{x_j}\Big|^2
 &\ds = |I_1|^2+2I_1I_2+|I_2|^2\ge 2I_1I_2.
 \ea
 \ee

By virtue of  \cite[the\ equation
(4.8)]{Fu-Yong-Zhang1},  a simple calculation
shows that

 \bel{1c1.3}
 \3n\1n\ba{ll}
 \displaystyle
 I_1I_2\3n&\ds=2\sum_{i,j=1}^m a^{ij}\ell_{x_i}v_{x_j}\(-\sum_{i,j=1}^m (a^{ij}v_{x_i})_{x_j}+\L v\)+\Psi v \(-\sum_{i,j=1}^m (a^{ij}v_{x_i})_{x_j}+\L v\)\\
 \ns&\ds
 = - \2n\sum_{j=1}^m\(2\sum_{i,i',j'=1}^ma^{ij} a^{i'j'}\ell_{x_{i'}}v_{x_i}v_{x_{j'}}
 -\sum_{i,i',j'=1}^ma^{ij}a^{i'j'}\ell_{x_i}v_{x_{i'}}v_{x_{j'}}-\L \sum_{i=1}^ma^{ij}\ell_{x_i} v^2\)_{x_j}\\
 \noalign{\ss}
 &\displaystyle
 \q+\sum_{i,j,i',j'=1}^m \left(2a^{ij'}\big(a^{i'j}\ell_{x_{i'}}\big)_{x_{j'}} -
 \big(a^{ij}a^{i'j'}\ell_{x_{i'}}\big)_{x_{j'}}\right)v_{x_i}v_{x_j}-
 \sum_{i,j=1}^m\big(\L a^{ij}\ell_{x_i}\big)_{x_j} v^2\\
 \ns&\ds\q -\sum_{i,j=1}^m
 \big(\Psi a^{ij}vv_{x_i}\big)_{x_j}+ \Psi \sum_{i,j=1}^m
 a^{ij}v_{x_i}v_{x_j}+\sum_{i,j=1}^m a^{ij}\Psi_{x_j}vv_{x_i}+\L \Psi v^2.
 \ea
 \ee
Combining   (\ref{1c1.3})  with (\ref{c1e5}), we
obtain the desired inequality (\ref{c1e2a}). \endpf

%%%%%%%%%%%%%%%%%%%%%%%%%%%%%%%%%%%%%%%%%%%%%%%%%%

\section{A Carleman estimate for  deterministic  hyperbolic equations in the $H^1$-norm}
\label{sec-H1}

%%%%%%%%%%%%%%%%%%%%%%%%%%%%%%%%%%%%%%%%%%%%%%%%%%%

This section is addressed to deriving a Carleman estimate for the following (deterministic) hyperbolic
equation:
\begin{eqnarray}{\label{system2}}
\left\{
\begin{array}{lll}\ds
\ds u_{tt} - \sum_{i,j=1}^n(b^{ij}u_{x_i})_{x_j}
= F  & {\mbox { in }} Q,
 \\
\ns\ds  u = 0 & \mbox{ on } \Si, \\[3mm]
\ns\ds  u(0) = u_0, \ u_{t}(0) = u_1 & \mbox{ in
} G,
\end{array}
\right.
\end{eqnarray}
where  $(u_0, u_1) \in H_0^1(G) \times L^2(G)$,
$F\in L^2(Q)$, and $b^{i j}$ $(i, j=1, \cdots,
n)$ satisfy  (\ref{gz1}), (\ref{gz2}) and the Condition \ref{condition of d}.

As in Theorem \ref{inobser}, we assume that  $T>T_0$ (defined in (\ref{fllz101})), and $G_0$, $\mu_0$ and $d(\cdot)$ are given in  (\ref{20160115e2})  and  Condition
\ref{condition of d}.
By (\ref{t0}), it is easy to see that  $
 \frac{T_0}{T}<\frac{\sqrt{\mu_0 M_0}}{3T}$. Hence, we can choose a constant $c_1\in (T_0/T, \min\{1, \frac{\sqrt{\mu_0 M_0}}{3T}\})$.
Now, for any given constant $c_0\in (0, 1)$ and parameter $\lambda>0$, we
choose the weight function $\th$ and the auxiliary
function $\Psi$ (appeared in Lemma \ref{c1t1}) as follows:
 \bel{gz4}\left\{\ba{ll}\ds
\th(t,x)=e^{\ell(t,x)},\qq \ell(t,x)=\l\phi(t,x),\qq \phi(t, x)=d(x)-c_1(t-T/2)^2,\\
\ns\ds \Psi(x)=\l\(\sum_{i,
j=1}^n(b^{ij}d_{x_i})_{x_j}-2c_1-c_0\).
\ea\right.
 \ee
We have  the
following global Carleman estimate for the
equation (\ref{system2}).

\bt\label{5.24-lm1} Assume that Condition
$\ref{condition of d}$  holds.   Then, there is a  positive constant
$\l_0$,  such that for any $T>T_0$ and $\l\geq \l_0$, any
solution $u$  to \eqref{system2} satisfies that
\begin{equation}\label{5.24-lm1-eq1}\ba{ll}\ds
\int_Q e^{2\lambda\phi }\big[\l(u_t^2 + |\nabla
u|^2) + \l^3 u^2\big] dxdt \leq C\[ \int_Q e^{2\lambda\phi }F^2
dxdt+\l^2\int_0^T\int_{G_0} e^{2\lambda\phi }(u_t^2 +\lambda^2
u^2)dxdt \]. \ea
\end{equation}
\et

\begin{remark} Notice that Theorem $\ref{5.24-lm1}$ is an improvement of \cite[Theorem 5.1]{Fu-Yong-Zhang1}.
 Indeed, in \cite[Theorem 5.1]{Fu-Yong-Zhang1}, the global Carleman estimate $(\ref{5.24-lm1-eq1})$ for  a deterministic hyperbolic equation
 was established under the additional condition that
 $u(0)=u(T)=0$ in $G$.    However,  this condition seems  too strong to be satisfied in applications $($see, for example
the  equation  $(7.5)$ in $\cite{Fu-Yong-Zhang1})$.  Therefore, it is necessary to establish
 the  global Carleman estimate  $(\ref{5.24-lm1-eq1})$
  without  this  restriction.\end{remark}
  
In the rest of this section, we  give a proof of Theorem
\ref{5.24-lm1}. 

\ms

\noindent{\bf Proof of Theorem \ref{5.24-lm1}.}
The proof is long and therefore we divide it into
four steps.

\ms

\noindent {\bf  Step 1.  A pointwise inequality
for hyperbolic operators. } In Lemma \ref{c1t1},
we choose $m=n+1$ and
 $\ds
 (a^{ij})_{m\times m}=\left(\ba{cc}
 -1&0\\0&(b^{ij})_{n\times n}
 \ea
 \right)
 $,  and $\th$, $\ell$, $\phi$ and $\Psi\equiv\Psi(x)$ being given as in (\ref{gz4}).  Then, by a simple calculation,  we have the following weighted inequality for the hyperbolic operator,  which is very similar to \cite[Corollary 4.2]{Fu-Yong-Zhang1}, except some different lower order terms.
 \bel{b3}
 \q\ba{ll}
\displaystyle
 e^{2\lambda\phi}\left|u_{tt}-\sum_{i,j=1}^n(b^{ij}u_{x_i})_{x_j}\right|^2+2\sum_{j=1}^n\Big(2\sum_{i, i', j'=1}^nb^{ij}
b^{i'j'}\ell_{x_{i'}}v_{x_i}v_{x_{j'}}-\sum_{i=1}^n b^{ij}\L \ell_{x_i} v^2
 \\
 \noalign{\ss}
 \displaystyle-\sum_{i, i', j'=1}^nb^{ij}b^{i'j'}\ell_{x_i}v_{x_{i'}}v_{x_{j'}}+\Psi v\sum_{i=1}^n
 b^{ij}v_{x_i}-2\ell_tv_t\sum_{i=1}^n b^{ij}v_{x_i}+\sum_{i=1}^nb^{ij}\ell_{x_i}v_t^2\Big)_{x_j}+2M_t\\
 \noalign{\ss}
 \displaystyle
 \ge2\Big(\ell_{tt}+\sum_{i, j=1}^n(b^{ij}\ell_{x_i})_{x_j}-\Psi\Big)v_t^2
 -8\sum_{i, j=1}^nb^{ij}\ell_{x_jt}v_{x_i}v_t +2\sum\limits_{i, j=1}^{n} b^{i j}\Psi_{x_j}v v_{x_i}+Bv^2\\
 \noalign{\ss}
 \displaystyle\q
 +2\sum_{i, j=1}^n \Big[b^{ij}\ell_{tt}
+\sum_{i',
j'=1}^n\(2b^{ij'}(b^{i'j}\ell_{x_{i'}})_{x_{j'}}
- (b^{ij}b^{i'j'}\ell_{x_{i'}})_{x_{j'}}\)+\Psi b^{ij}
 \Big]v_{x_i}v_{x_j},   \ea
 \ee
where
\bel{b7}
\left\{
 \ba{ll}
\ds \L =\ell _t^2-\ell _{tt}-\sum_{i, j=1}^n
 (b^{ij}\ell_{x_i}\ell_{x_j}-b^{ij}_{x_j}\ell_{x_i}
 -b^{ij}\ell_{x_ix_j})-\Psi,\\
 \ns\ds M=\ell_t\(v_t^2+\sum_{i, j=1}^n
 b^{ij}v_{x_i}v_{x_j}\)-2\sum_{i, j=1}^nb^{ij}\ell_{x_i}v_{x_j}v_t
 -\Psi vv_t+\L \ell_tv^2,\\
\ns\ds
\ds B=2\Big(\L \Psi+(\L \ell_t)_t-
 \sum_{i, j=1}^n \big(\L b^{ij}\ell_{x_i}\big)_{x_j}\Big).
 \ea
\right.
\ee

\ms

\noindent {\bf  Step 2.  Estimates on  ``the
energy terms". }  First, by  the definitions of
$\Psi$ and $\ell$,   it is easy to show that
  \bel{gz3}\ba{ll}\ds
  2\Big(\ell_{tt}+\sum_{i, j=1}^n(b^{ij}\ell_{x_i})_{x_j}-\Psi\Big)=-4\lambda c_1+2\sum\limits_{i, j=1}^{n} (\lambda b^{i j}d_{x_i})_{x_j}-2\Psi=2\lambda c_0.
 \ea
 \ee
Further,
\begin{eqnarray}\label{fllz1}
\begin{array}{rl}
&2\displaystyle\sum\limits_{i, j=1}^n
\Big[b^{ij}\ell_{tt} +\sum\limits_{i',
j'=1}^n\(2b^{ij'}(b^{i'j}\ell_{x_{i'}})_{x_{j'}}
-
 (b^{ij}b^{i'j'}\ell_{x_{i'}})_{x_{j'}}\)+\Psi b^{ij}
 \Big]v_{x_i}v_{x_j}\\[3mm]
&=2\displaystyle\sum\limits_{i, j=1}^n
\Big(-2\lambda c_1 b^{ij}
+\displaystyle\sum\limits_{i', j'=1}^n  \lambda b^{ij}(b^{i' j'}d_{x_{i'}})_{x_{j'}}-2\lambda c_1b^{i j}-c_0\lambda b^{i j}\\[3mm]
&\quad\quad\quad\quad+2\lambda
\displaystyle\sum\limits_{i', j'=1}^{n} b^{i
j'}(b^{i' j}d_{x_{i'}})_{x_{j'}}-\lambda
\sum\limits_{i', j'=1}^{n} (b^{i j}b^{i'
j'}d_{x_{i'}})_{x_{j'}}
 \Big)v_{x_i}v_{x_j}\\[3mm]
&\geq 2\lambda\mu_0\displaystyle\sum\limits_{i,
j=1}^{n} b^{i
j}v_{x_i}v_{x_j}-(8c_1+2c_0)\lambda\sum\limits_{i,
j=1}^{n} b^{i j}v_{x_i}v_{x_j}
=2\lambda(\mu_0-4c_1-c_0)\sum\limits_{i,
j=1}^{n} b^{i j}v_{x_i}v_{x_j}.
\end{array}
\end{eqnarray}
Further, by (\ref{gz4}) and (\ref{b7}), we obtain
that
 \begin{eqnarray*}\ba{ll}\ds
 &\L\ds=\ell _t^2-\ell _{tt}-\sum_{i, j=1}^n
 (b^{ij}\ell_{x_i}\ell_{x_j}-b^{ij}_{x_j}\ell_{x_i}
 -b^{ij}\ell_{x_ix_j})-\Psi\\
 \ns&=\l^2\Big[c_1^2(2t-T)^2-\displaystyle\sum\limits_{i,j=1}^nb^{ij}d_{x_i}d_{x_j}\Big]+O(\l).
 \ea\end{eqnarray*}
Hence,
 \bel{c4}\ba{ll}\ds
B=2\l^3\[(4c_1+c_0)\sum_{i, j=1}^nb^{ij}d_{x_i}d_{x_j}+\sum_{i, j=1}^nb^{ij}d_{x_i}\(\sum_{i', j'=1}^nb^{i'j'}d_{x_{i'}}d_{x_{j'}}\)_{x_j}\\[3mm]
\ds\quad\quad\quad\quad-(8c_1+c_0)c^2_1(2t-T)^2\]+O(\l^2).
\ea \ee Proceeding  the same analysis  as
(11.6)-(11.8) in  \cite{Fu-Yong-Zhang1}, we
have that
 \bel{otc2}
 \sum_{i,j=1}^n\sum_{i',j'=1}^nb^{ij}d_{x_i}(b^{i'j'}d_{x_{i'}}d_{x_{j'}})_{x_j}\ge \mu_0\sum_{i,j=1}^nb^{ij}d_{x_i}d_{x_j}.
 \ee
By (\ref{c4}), (\ref{otc2}) and
(\ref{t0}), noticing that  $c_1<\frac{\sqrt{\mu_0 M_0}}{3T}$, we find that for any $T>T_0$,
 \bel{pb8}
B\ge
2\l^3(4c_1+c_0)\sum_{i,j=1}^nb^{ij}d_{x_i}d_{x_j}+O(\l^2).
 \ee
 
On the other hand, by (\ref{t0}), it is easy to see that
 $$
  \frac{\mu_0M_0}{ {8c_1+c_0}}> \frac{\mu_0 M_0}{9}>T_0^2\ge4M_1\ge 4M_0,
 $$
which implies that
 \bel{1221}
 \mu_0-4c_1-c_0>\mu_0-32c_1-4c_0>0.
 \ee
Therefore, combining (\ref{gz3}),  (\ref{fllz1})
(\ref{pb8}) and (\ref{1221}) with (\ref{b3}), we conclude
that for any $T>T_0$, there is a $\l_1>0$ and $c^*>0$, such that
for any $\l\ge\l_1$,
  \bel{gz5}\ba{ll}\ds
  2\Big(\ell_{tt}+\sum_{i, j=1}^n(b^{ij}\ell_{x_i})_{x_j}-\Psi\Big)v_t^2
 -8\sum_{i, j=1}^nb^{ij}\ell_{x_jt}v_{x_i}v_t +2\sum\limits_{i, j=1}^{n} b^{i j}\Psi_{x_j}v v_{x_i}+Bv^2\\
 \noalign{\ss}
 \displaystyle\q
 +2\sum_{i, j=1}^n \Big[b^{ij}\ell_{tt}
+\sum_{i',
j'=1}^n\(2b^{ij'}(b^{i'j}\ell_{x_{i'}})_{x_{j'}}
-
 (b^{ij}b^{i'j'}\ell_{x_{i'}})_{x_{j'}}\)+\Psi b^{ij}
 \Big]v_{x_i}v_{x_j}\\[5mm]
 \ge c^*\l (v_t^2+|\nabla v|^2+\l^2v^2).
 \ea
 \ee%
Integrating (\ref{gz5}) in $Q$ and  noting that
$v=e^{\lambda\phi} u$  on $\Si$,  by  (\ref{system2}) and
(\ref{b3}), we obtain  that
  \bel{gz6}\ba{ll}\ds
  c^*\l\int_Q (v_t^2+|\nabla v|^2+\l^2v^2)dxdt\\
  \ns\ds \le \int_Q e^{2\lambda\phi}|F|^2dxdt+2\int_Q M_tdxdt
  +2\l d_0\int_{\Si_0}\sum_{i,j=1}^nb^{ij}\n^i\n^je^{2\lambda\phi}\|\frac{\pa u}{\pa\nu}\|^2d\Si.
  \ea\ee
Here  we  use the following identity:
  $$\ba{ll}\ds
\int_{\Si} \sum_{i, j, i', j'=1}^n\(2b^{ij}
b^{i'j'}\ell_{x_{i'}}v_{x_i}v_{x_{j'}}
 -b^{ij}b^{i'j'}\ell_{x_i}v_{x_{i'}}v_{x_{j'}}\)\cdot\nu^jd\Si\\
 \ns\ds=\l\int_{\Si}\sum_{i,j=1}^nb^{ij}\n^i\n^j\sum_{i',j'=1}^nb^{i'j'}d_{x_{i'}}\n^{j'}\|\frac{\pa v}{\pa\nu}\|^2d\Si.
 \ea$$
\noindent {\bf  Step 3.  Estimates on ``the
spatial  boundary term". } Let us estimate the
last term of (\ref{gz6}).   Similar to the proof of (11.15)
in \cite{Fu-Yong-Zhang1}, we choose  functions
$h_0\in C^1(\cl{G}; [0,1]^n)$  and $\rho\in
C^2(\cl{G};[0,1])$,  such that $\ds h_0=\nu$ on
$\G$, and for the same $\delta$  appeared in (\ref{20160115e2}),
\begin{eqnarray*}
 \left\{\ba{ll}
 \rho(x)\equiv 1, \q x\in \cO_{\d/3}(\G_0)\cap G, \\[3mm]
 \rho(x)\equiv 0, \q x\in G\setminus\cO_{\d/2}(\G_0).
 \ea\right.
 \end{eqnarray*}
Let $h=h_0\rho e^{2\lambda\phi}$.  Then,  by  \cite[Lemma
3.2]{Fu-Yong-Zhang1} (with $g$ and $a^{ij}$
replaced by $h$ and $b^{ij}$, respectively),
 \begin{eqnarray*}
 \qq\ba{ll}
  \displaystyle
 \int_\Si\sum_{i, j=1}^nb^{ij}\nu^i\nu^j\rho  e^{2\lambda\phi}\|\frac{\pa
u}{\pa\nu}\|^2 dxdt\\
 \ns\ds=- \int_Q\Big[2\(F h\cdot\nabla u-(u_th\cdot\nabla
 u)_t+u_th_t\cdot\nabla u-(\nabla\cdot
 h)u_t^2\)\\
 \noalign{\ss}
 \displaystyle\qq\qq\qq
 -\sum_{i, j,k=1}^nb^{ij}u_{x_i}u_{x_k}\frac{\pa
 h^k}{\pa x_j}+\sum_{i, j=1}^nu_{x_i}u_{x_j}\nabla\cdot(b^{ij}h)\Big]dxdt.
 \ea
 \end{eqnarray*}
 Therefore, it is easy to check that

  \begin{eqnarray*}
 \qq\ba {ll}
  \displaystyle
 \int_\Si\sum_{i, j=1}^nb^{ij}\nu^i\nu^j\rho e^{2\lambda\phi}\|\frac{\pa
u}{\pa\nu}\|^2 dxdt\\
\ns\ds\leq \frac{C}{\lambda}\int_Q e^{2\lambda\phi}F^2dxdt+\lambda \ds\int^T_0\int_{\cO_{\d/2}(\G_0)\cap G} e^{2\lambda\phi} |\nabla u|^2dxdt+2\ds\int_G u_t h\cdot\nabla udx\Big|^T_0\\
\ns\ds\quad
+4\sqrt{n} \lambda\int^T_0\int_{\cO_{\d/2}(\G_0)\cap G} e^{2\lambda\phi} c_1T|u_t||\nabla u|dxdt+C\lambda\int^T_0\int_{\cO_{\d/2}(\G_0)\cap G} e^{2\lambda\phi} u_t^2dxdt\\
\ns\quad\ds+\(6\sqrt{n}\lambda|\nabla d|\sum\limits_{i, j=1}^{n}|b^{i j}|_{C(\overline{G})}+C\)\int^T_0\int_{\cO_{\d/2}(\G_0)\cap G}e^{2\lambda\phi}|\nabla u|^2dxdt\\
\ns\ds\leq C\(\frac{1}{\lambda}\int_Q e^{2\lambda\phi}F^2dxdt+\lambda\int^T_0\int_{\cO_{\d/2}(\G_0)\cap G} e^{2\lambda\phi} u_t^2dxdt\)+2\ds\int_G u_t h\cdot\nabla udx\Big|^T_0\\
\ns\ds\quad+6\sqrt{n}\lambda\(|\nabla d|\sum\limits_{i, j=1}^{n}|b^{i j}|_{C(\overline{G})}+1\)\int^T_0\int_{\cO_{\d/2}(\G_0)\cap G}e^{2\lambda\phi}|\nabla u|^2dxdt.
 \ea
 \end{eqnarray*}
By  (\ref{gz6}) and the above inequality,  we get that
  \bel{fllz2}\ba{ll}\ds
  c^*\l\int_Q (v_t^2+|\nabla v|^2+\l^2v^2)dxdt\\
  \ns\ds\leq C\(\int_Q e^{2\lambda\phi}F^2dxdt+\lambda^2 \int^T_0\int_{G_0} e^{2\lambda\phi} u^2_t dxdt\)\\
  \ns\ds\q+4\lambda d_0 \int_G u_t h\cdot\nabla udx\Big|^T_0+2\int_G \tilde{M} dx\Big|^T_0\\
  \ns\ds\quad+
12\sqrt{n}\lambda^2 d_0\(|\nabla d|\sum\limits_{i, j=1}^{n}|b^{i j}|_{C(\overline{G})}+1\)\int^T_0\int_{\cO_{\d/2}(\G_0)\cap G}e^{2\lambda\phi}|\nabla u|^2dxdt,
  \ea
  \ee
 where
 $\tilde{M}=M+2\lambda d_0 u_t h\cdot\nabla u$.

\medskip

Next, let us estimate the last term in
(\ref{fllz2}). Put $
 \eta(t, x)=\rho_1^2 e^{2\lambda\phi},
$
 where $\rho_1\in C^2(\cl{G};[0,1])$ satisfies that
\begin{eqnarray*}
 \left\{\ba{ll}
% 0\leq \rho_1(x)\leq  1, \q x\in\oO, \\
\rho_1(x)\equiv 1, \q x\in \cO_{\d/2}(\G_0)\cap G, \\[3mm]
 \rho_1(x)\equiv 0, \q x\in G\setminus G_0.
 \ea\right.
 \end{eqnarray*}
By (\ref{system2}), we have that
 \begin{eqnarray*}\;\;\;\;
 \ba{ll}
 &\displaystyle
\int_Q \eta u Fdxdt=\int_Q \eta u\(u_{tt}-\sum_{i,j=1}^n(b^{ij}u_{x_i})_{x_j}\)dxdt\\
 \noalign{\ss}
 &\displaystyle
 =\int_Q (\eta uu_t)_tdxdt-\int_Qu_t(\eta_t u+\eta u_t)dxdt
\\
\ns&\ds\q+\int_Q\eta
\sum_{i,j=1}^nb^{ij}u_{x_i}u_{x_j}dxdt+\int_Q
u\sum_{i,j=1}^nb^{ij}u_{x_i}\eta_{x_j}dxdt.
 \ea
 \end{eqnarray*}
 This implies that
 \bel{d35}\ba{ll}\ds
\int_0^T\int_{\cO_{\d/2}(\G_0)\cap G} e^{2\lambda\phi}|\nabla u|^2dxdt\\
\ns\ds
 \leq C\[\frac{1}{\l^2}\int_Q  e^{2\lambda\phi}|F|^2dxdt+\int_0^T\int_{G_0} e^{2\lambda\phi}(\l^2u^2+u_t^2)dxdt\]-\ds\frac{1}{s_0}\int_Q (\eta uu_t)_tdxdt.
  \ea\ee
 Combining (\ref{fllz2}) with (\ref{d35}), we end up with
  \bel{otc6}\ba{ll}\ds
 c^* \l\int_Q (v_t^2+|\nabla v|^2+\l^2v^2)dxdt\\
  \ns\ds \le C\int_Q e^{2\lambda\phi}|F|^2dxdt+\int_G {\overline M}dx\Big|^T_0+C\l^2\int_0^T\int_{G_0} e^{2\lambda\phi}(u_t^2 +\lambda^2
 u^2)dxdt,
  \ea\ee
where
 \begin{eqnarray*}\ba{ll}\ds
 \overline M=M+2d_0\l u_t h\cdot\nabla u-\frac{1}{s_0}\[12\sqrt{n}\lambda^2 d_0\(|\nabla d|\sum\limits_{i, j=1}^{n}|b^{i j}|_{C(\overline{G})}+1\)\]\eta uu_t.
 \ea\end{eqnarray*}

 \smallskip

\noindent {\bf  Step 4.  Estimates on ``the time
boundary term". } Let us  estimate  $\overline
M(0,x)$ and $\overline M(T,x)$, respectively. By
 (\ref{t0}) and the
definition of $M$ in   (\ref{b7}), we have that
  \begin{eqnarray*}\ba{rl}\ds
   M(0,x)
  &\ds\ge\[\l c_1T\(v_t^2+\sum_{i, j=1}^n
 b^{ij}v_{x_i}v_{x_j}\)-2\l\sum_{i, j=1}^nb^{ij}d_{x_i}v_{x_j}v_t\]\|_{t=0}\\
 \ns&\ds\q+\[O(\l^2)v^2-v_t^2+c_1T\l^3\(c_1^2T^2-\sum_{i,j=1}^nb^{ij}d_{x_i}d_{x_j}\)v^2\]\|_{t=0}\\
 \ns&\ds\ge \l \left(c_1 T-\(\sum\limits_{i, j=1}^{n} b^{i j}d_{x_i}d_{x_j}\)^{\frac{1}{2}} \right)\(v_t^2+\sum_{i,j=1}^nb^{ij}v_{x_i}v_{x_j}\)\|_{t=0}\\
 \ns&\ds\q+\[O(\l^2)v^2-v_t^2+c_1T\l^3\(c_1^2T^2-\sum_{i,j=1}^nb^{ij}d_{x_i}d_{x_j}\)v^2\]\|_{t=0}.
 \ea\end{eqnarray*}
 Noting that by  (\ref{fllz101})  and $c_1>T_0/T$, we have
  \bel{0lfz0}
  c_1T>2\sqrt{M_1}\ge 2\(\sum_{i,j=1}^n
b^{ij}(x)d_{x_i}(x)d_{x_j}(x)\)^{\frac{1}{2}}\ge 2\sqrt{s_0}|\nabla d|.
  \ee
 This implies that
  \bel{lfz0}
M(0,x)\ge \[\frac{1}{2}\l c_1T\min\{1,s_0\}(v_t^2+|\nabla v|^2)+\frac{3}{4}\l^3c_1^3T^3v^2+O(\l^2)v^2-v_t^2\]\|_{t=0}.
  \ee
Further,
  \bel{lfz1}\ba{ll}\ds
 2d_0\l u_t h\cdot\nabla u\|_{t=0}\ge - \sqrt{n}d_0\l e^{2\lambda\phi}(u_t^2+|\nabla u|^2)\|_{t=0}\\
 \ns\ds\ge - 2\sqrt{n}d_0\l \Big( v_t^2+ \l^2c_1^2T^2v^2+|\nabla v|^2+\lambda^2|\nabla d|^2 v^2\Big)\|_{t=0}.
 \ea
 \ee

 On the other hand,
  \begin{eqnarray}\label{10*}\ba{ll}\ds
  -\frac{1}{s_0}\[12\sqrt{n}\lambda^2 d_0\(|\nabla d|\sum\limits_{i, j=1}^{n}|b^{i j}|_{C(\overline{G})}+1\)\]\eta uu_t\\
  \ns\ds=-\frac{1}{s_0}\[12\sqrt{n}\lambda^2 d_0\(|\nabla d|\sum\limits_{i, j=1}^{n}|b^{i j}|_{C(\overline{G})}+1\)\]\rho_1^2 \(v v_t-\lambda c_1Tv^2-\frac{v_t^2}{\l c_1T}+\frac{v_t^2}{\l c_1T}\)\Big|_{t=0}\\
    \ns\ds\geq  -\frac{1}{s_0 c_1T}\lambda\[12\sqrt{n} d_0\(|\nabla d|\sum\limits_{i, j=1}^{n}|b^{i j}|_{C(\overline{G})}+1\)\]v_t^2\Big|_{t=0}.
  \ea\end{eqnarray}
 Therefore, by (\ref{lfz0})-(\ref{10*}), we get that
  \bel{lfz3}\ba{ll}\ds
  \overline M(0,x)\ge \(\l F_1v_t^2+\l F_2|\nabla v|^2+\l^3F_3v^2+O(\l^2)v^2+O(1)v_t^2\)\|_{t=0},
  \ea\ee
 where
  \bel{ll}\left\{\ba{ll}\ds
 F_1=\frac{1}{2} c_1T\min\{1,s_0\}-2\sqrt{n}d_0-\frac{1}{s_0 c_1T} \[12\sqrt{n} d_0\(|\nabla d|\sum\limits_{i, j=1}^{n}|b^{i j}|_{C(\overline{G})}+1\)\],\\
 \ns\ds F_2=\frac{1}{ 2} c_1T\min\{1,s_0\}-2\sqrt{n}d_0,\\
 \ns\ds F_3=\frac{3}{4}c_1^3T^3-2\sqrt{n}d_0(c_1^2T^2+|\nabla d|^2).
 \ea\right.\ee
By (\ref{0lfz0}) and  (\ref{fllz101}),  for any $T>T_0$,  it holds that $c_1T>1$,
 $$
 F_1\ge \frac{1}{ 2}c_1T\min\{1,s_0\}-2\sqrt{n}d_0-\frac{6\sqrt{n}d_0}{s_0^{3/2}} \sum\limits_{i, j=1}^{n}|b^{i j}|_{C(\overline{G})}
 -\frac{12\sqrt{n}d_0}{s_0}>0,
 $$
 and therefore, $F_2>0$. Moreover,
  $$
 F_3\ge \frac{3}{4}c_1^3T^3-2\sqrt{n}d_0c_1^2T^2\(1+\frac{1}{ 4s_0}\)=\frac{3}{4}c_1^2T^2\[c_1T-\frac{8}{ 3}\sqrt{n}d_0\(1+\frac{1}{ 4s_0}\)\]>0,
  $$
 where we use the following fact:
  $$\ba{ll}\ds
 &\ds c_1T>\frac{4\sqrt{n}}{ \min\{1,s_0\}}d_0\ds=\frac{8\sqrt{n}}{ 3\min\{1,s_0\}}d_0+\frac{4\sqrt{n}}{ 3\min\{1,s_0\}}d_0\\
  \ns&\ds \ge \frac{8\sqrt{n}}{ 3}d_0+\frac{4\sqrt{n}}{ 3 s_0 }d_0= \frac{8}{ 3}\sqrt{n}d_0\(1+\frac{1}{ 2s_0}\).
  \ea$$
 Finally, by (\ref{fllz101}),  one can find constants $C>0$ and $\l_2>0$ such that for any $\l\ge\l_2$,
 \bel{gz9}\ba{ll}\ds
  \overline M(0,x)\ge 0.
 \ea\ee
Meanwhile, noting that $\ell(T,x)=-\ell(0,x)$,
we have that  there is a constant $\l_3>0$ such
that for any $\l\ge\l_3$,
  \bel{gz10}\ba{ll}\ds
  \overline M(T,x)\le 0.
 \ea\ee
 Combining  (\ref{gz9}) and (\ref{gz10}) with (\ref{otc6}),  and noting that $v=e^{\lambda\phi} u$,  for any $\l\ge\l_0=\max\{\l_1,\l_2,\l_3\}$, we end up with the desired estimate
(\ref{5.24-lm1-eq1}). This completes the proof of Theorem
\ref{5.24-lm1}.\endpf

\ms

\section{An auxiliary optimal control problem}\label{ss6}

In this section, as a preliminary, we analyze an auxiliary optimal
control problem.
Some ideas are taken from \cite[pp. 190--199]{I}
and \cite[Proposition 6.1]{Fu-Yong-Zhang1}.

\smallskip

Let $y\in L^2_\mathbb{F}(\O;C([0,T];L^2(G)))$
satisfy $y(0)=y(T)=0$ in $G$, $\dbP$-a.s. For
any $K>1$, let $\varrho\equiv\varrho^K(x)\in
C^2(\cl{G})$, such that
$\ds\min_{x\in\cl{G}}\varrho(x)=1$ and
\begin{equation}\label{csss}
\varrho(x)= \left\{\begin{array}{ll}
\ds 1 &\mbox{ if $x\in G_0$}, \\
\ns\ds   K&\mbox{ if ${\rm dist}(x, G_0)\geq
\ds\frac{1}{{\rm ln} K}$}.
\end{array}
\right.
\end{equation}
For any integer $m\geq3$, let $h=\ds\frac{T}{m}$, and set
\begin{equation}\label{csum}
y_m^j\equiv y_m^j(x)=y(jh,x),\q\q\phi_m^j\equiv \phi_m^j(x)=\phi(jh,x),\q
j=0,1,\cds,m,
\end{equation}
where $\phi(t, x)=d(x)-c_1(t-T/2)^2$ (see (\ref{gz4})).
Consider the following system:
\begin{equation}\label{5.23-eq3}
\left\{
\begin{array}{ll}
\ds
\mE\(\frac{z_m^{j+1}-2z_m^j+z_m^{j-1}}{h^2}\ \!\Big|\ \!\cF_{jh}\)-
\sum_{j_1,j_2=1}^n\pa_{x_{j_2}}(b^{j_1j_2}\pa_{x_{j_1}}z_m^j)\\
\ns\ds\quad\quad=
\mE\(\frac{r_{1m}^{j+1}\!-\!r_{1m}^j}{h}\ \!\Big|\ \!\cF_{jh}\)\!+\!r_{2m}^j
\!+\!\l y_m^je^{2\l \phi_m^j}+r_m^j\ \
(1\le j\le m-1) &\mbox{in }\ G, \\[3mm]
\ns \ds
z_m^j=0\q (0\le j\le m)  &\mbox{on }\ \G, \\[3mm]
\ns\ds z_m^0=z_m^m=r_{2m}^0=r_{2m}^m=r_m^0=r_m^m=0, \q r_{1m}^0=r_{1m}^1 &\mbox{in }\ G.
\end{array}
\right.
\end{equation}
Here $(r_{1m}^j, r_{2m}^j, r_m^j)\in
\Big(L^2_{\cF_{jh}}(\O;L^2(G))\Big)^3$
($j=0,1,\cdots,m$) are  controls. The set of
admissible sequences for (\ref{5.23-eq3}) is
defined by
$$
\3n\begin{array}{ll}\ds
\cA_{ad}\ds\!=\!\Big\{\{(z_m^j,r_{1m}^j,
r_{2m}^j, r_m^j)\}_{j=0}^m; \,(z_m^j,r_{1m}^j,
r_{2m}^j, r_m^j) \in
L^2_{\cF_{jh}}(\O;H_0^1(G))\times
\Big(L^2_{\cF_{jh}}(\O;L^2(G))\Big)^3
\\
\ns\ds  \hspace{5.7cm} \mbox{ and }
\{(z_m^j,r_{1m}^j, r_{2m}^j, r_m^j)\}_{j=0}^m
\mbox{
 solves
} (\ref{5.23-eq3})\Big\}.
\end{array}
$$
Since $\{(0, 0, 0, -\l y_m^je^{2\l
\phi_m^j})\}_{j=0}^m\in\cA_{ad}$, it follows that
$\cA_{ad}\neq\emptyset$.

\medskip

Next,  define a cost functional as follows:
\begin{equation}\label{5.23-eq0}
\begin{array}{ll} \ds
J(\{(z_m^j,r_{1m}^j, r_{2m}^j, r_m^j)\}_{j=0}^m)
\3n&\ds =\frac{h}{2}\mE\int_G \varrho
\frac{|r_{1m}^m|^2}{\l^2}e^{-2\l
\phi_m^m}dx+\frac{h}{2}\mE\sum_{j=1}^{m-1}\[\int_G|z_m^j|^2
e^{-2\l \phi_m^j}dx \\
\ns&\ds \q +\int_G \varrho
\(\frac{|r_{1m}^j|^2}{\l^2}+\frac{|r_{2m}^j|^2}{\l^4}\)e^{-2\l
\phi_m^j}dx +K\int_G|r_m^j|^2 dx\],
\end{array}
\end{equation}
and consider the following optimal control
problem: Find a $\{(\hat z_m^j,\hat r_{1m}^j,
\hat r_{2m}^j, \hat r_m^j)\}_{j=0}^m\in
\cA_{ad}$,  such that
\begin{equation}\label{5.23-eq1}
J(\{(\hat z_m^j,\hat r_{1m}^j, \hat r_{2m}^j,
\hat r_m^j)\}_{j=0}^m)  =
\min_{\{(z_m^j,r_{1m}^j, r_{2m}^j,
r_m^j)\}_{j=0}^m\in \cA_{ad}}
J(\{(z_m^j,r_{1m}^j, r_{2m}^j,
r_m^j)\}_{j=0}^m).
\end{equation}
Notice that for any $\{(z_m^j,r_{1m}^j, r_{2m}^j,
r_m^j)\}_{j=0}^m\in\cA_{ad}$, by the standard
regularity results for elliptic equations, we
have that $z_m^j\in L^2_{\cF_T}(\O;H^2(G)\cap
H_0^1(G))$.

We have the following result.

\begin{proposition}\label{5.24-prop1}
For any  $K>1$ and $m\ge3$,  the problem
$(\ref{5.23-eq1})$ admits a unique solution
$\{(\hat z_m^j,\hat r_{1m}^j, \hat r_{2m}^j,
\hat r_m^j)\}_{j=0}^m\in \cA_{ad}$, $($which
depends on $K)$.  Furthermore, define
\begin{equation}\label{csxg1}
p_m^j\equiv p_m^j(x)\=K \hat r_m^j(x),\qq 0\le
j\le m.
\end{equation}
Then, 
\begin{equation}\label{csxgg2}
\left\{
\begin{array}{ll}\ds
\hat z_m^0=\hat z_m^m=p_m^0=p_m^m=0 \mbox{ in
}\ G,\\[3mm]
\ns\ds \hat z_m^j, \  p_m^j \in
L^2_{\cF_{jh}}(\O;H^2(G)\cap H_0^1(G)),\quad
1\le j\le m-1.
\end{array}
\right.
\end{equation}
Also,  the following optimality conditions hold:
\begin{equation}\label{5.23-eq4}
\left\{
\begin{array}{ll}
\ds\frac{p_m^j-p_m^{j-1}}{ h}+ \varrho
\frac{\hat r_{1m}^j}{\l^2}e^{-2\l \phi_m^j}=0 &
\mbox{ in
}G,\\
\ns\ds p_m^j-\varrho \frac{\hat
r_{2m}^j}{\l^4}e^{-2\l \phi_m^j}=0 &
\mbox{ in }G,\\
\end{array}
\qq\qq1\le j\le m,
\right.
\end{equation}
\begin{equation}\label{cs112}
\!\left\{
\!\begin{array}{ll}
\ds
\mE\(\frac{p_m^{j+1}-2 p_m^j +p_m^{j-1}}{h^2}\ \!\Big|\ \!\cF_{jh}\)\\
\ns\ds\qq
-\sum_{j_1,j_2=1}^n\pa_{x_{j_2}}\big(b^{j_1j_2}\pa_{x_{j_1}}
p_m^j \big)
+e^{-2\l \phi_m^j} \hat z_m^j=0& \mbox{ in }G,\\[3mm]
\ds p_m^j=0 &\mbox{ on }\ \G,
\end{array}
\quad 1\le j\le m-1. \right.
\end{equation}
Moreover, there is a constant $C=C(K,\l)>0$,
independent of $m$, such that
\begin{equation}\label{5.23-eq6}
h\mE\sum_{j=1}^{m-1}\int_G\[|\hat z_m^j|^2+|\hat
r_{1m}^j|^2+|\hat r_{2m}^j|^2+K|\hat
r_m^j|^2\]dx+h\mE\int_G |\hat r_{1m}^m|^2dx\le
C,
\end{equation}
and
\begin{equation}\label{5.24-eq7}
\begin{array}{ll}\ds
h\mE\sum_{j=0}^{m-1} \int_G\Big\{\frac{\big(\hat
z_m^{j+1}-\hat
z_m^j\big)^2}{h^2}+\frac{[\mE(\hat
r_{1m}^{j+1}-\hat
r_{1m}^j\ \!|\ \!\cF_{jh})]^2}{h^2}\\
\ns\ds \qq\qq\q +\frac{(\hat r_{2m}^{j+1}-\hat
r_{2m}^j)^2}{h^2} +K\frac{\big(\hat r_m^{j+1} -
\hat r_m^j\big)^2}{h^2}\Big\}dx\le C.
\end{array}
\end{equation}
\end{proposition}

We refer to  Appendix A for a proof of this
proposition.

\section{Global Carleman estimate for stochastic hyperbolic equations  in the
$L^2$-space}\label{ss7}

We define a formal differential operator $\cA$ by
 $$
 \cA\=\frac{\pa^2}{\pa{t^2}}-\sum_{i,j=1}^n\frac{\pa}{\pa{x_j}}\(b^{ij}\frac{\pa}{\pa{x_i}}\).
 $$

In order to prove Theorem \ref{inobser}, we need
the following global Carleman estimate for
stochastic hyperbolic equations.
\begin{theorem}\label{5.24-th1}
Assume that the Condition $\ref{condition of
d}$ holds. Let $T_0$  be given by (\ref{fllz101}).  Then there exists a
$\lambda_0^*>0$ such that for any $T>T_0$, $\l\ge\l_0^*$,
and any $y\in L^2_\dbF(\O;C([0,T];L^2(G)))$
satisfying $y(0)=y(T)=0$ in $G$ and
\begin{equation}\label{cspb6}
\ba{ll}
\ds \mE\big(y, \cA\eta\big)_{L^2(Q)}\\\ns
\ds =\mE\langle
b_1 y + f, \eta\rangle_{H^{-1}(Q),H_0^1(Q)}, \q\forall\;\eta\in
L^2_\dbF(\Omega; H^1_0(Q))\mbox{ with }
\cA\eta\in
L^2_\mathbb{F}(0, T; L^2(G)),
\ea
\ee
it holds that
\begin{equation}\label{5.24-eq9}
\begin{array}{ll}\ds
&\ds\l\mE\int_Q  e^{2\lambda\phi} y^2  dxdt\\
\ns&\ds \leq C\(\mE
|e^{\lambda\phi} f|_{H^{-1}(Q)}^2 +|e^{\lambda\phi} b_1
y|^2_{L^2_\dbF(0,T;H^{-1}(G))}  +\l^2
\mE\int_0^T\int_{G_0} e^{2\lambda\phi} y^2  dxdt\).
\end{array}
\end{equation}
\end{theorem}
\noindent {\bf Proof of Theorem \ref{5.24-th1}.
} We borrow some idea from \cite{Fu-Yong-Zhang1, I, liu0}. The whole  proof is divided into six steps.

\medskip

\noindent {\bf Step 1.} First, recall the
functions $\{(\hat z_m^j,\hat r_{1m}^j, \hat
r_{2m}^j, \hat r_m^j)\}_{j=0}^m$ in Proposition
\ref{5.24-prop1}.  For $m=2^i$ $(i=2, 3,
\cdots)$, we define
\begin{equation}\label{5.24-eq2}
\begin{array}{ll}\ds
\tilde z^m(t,x) =\frac{1}{
h}\sum_{j=0}^{m-1}\mE\(\Big\{(t-jh)\hat
z_m^{j+1}(x)-\big[t-(j+1)h\big]\hat
z_m^j(x)\Big\}\ \!\Big|\ \!\cF_{jh}\)\chi_{(jh,(j+1)h]}(t),\\
\ns\ds\tilde r_1^m(t,x) =\frac{1}{
h}\sum_{j=0}^{m-1}\mE\(\Big\{(t-jh)\hat
r_{1m}^{j+1}(x)-\big[t-(j+1)h\big]\hat
r_{1m}^j(x)\Big\}\ \!\Big|\ \!\cF_{jh}\)\chi_{(jh,(j+1)h]}(t),\\
 \ns\ds
\tilde r_2^m(t,x) =\frac{1}{
h}\sum_{j=0}^{m-1}\mE\(\Big\{(t-jh)\hat
r_{2m}^{j+1}(x)-\big[t-(j+1)h\big]\hat
r_{2m}^j(x)\Big\}\ \!\Big|\
\!\cF_{jh}\)\chi_{(jh,(j+1)h]}(t), \\
\ns\ds
\tilde r^m(t,x) =\frac{1}{
h}\sum_{j=0}^{m-1}\mE\(\Big\{(t-jh)\hat
r_{m}^{j+1}(x)-\big[t-(j+1)h\big]\hat
r_{m}^j(x)\Big\}\ \!\Big|\
\!\cF_{jh}\)\chi_{(jh,(j+1)h]}(t).
\end{array}
\end{equation}
By (\ref{5.23-eq6}) and (\ref{5.24-eq7}), there
is
 a subsequence  of
$\big\{\big(\tilde z^m, \tilde r_1^m, \tilde
r_2^m, \tilde r^m\big)\big\}_{m=2}^\infty$
(still denoted by itself), such that for some
$(\tilde z, \tilde r_1, \tilde r_2, \tilde r)
\in \big(L^2_{\dbF}(\O;H^1(0,T;L^2(G)))\big)^4$,
as $m\to\infty$,
\begin{equation}\label{5.26-eq11}
\big(\tilde z^m, \tilde r_1^m, \tilde r_2^m,
\tilde r^m\big)\rightarrow (\tilde z, \tilde
r_1, \tilde r_2, \tilde r) \mbox{ weakly in }
 \big(L^2_{\dbF}(\O;H^1(0,T;L^2(G)))\big)^4.
\end{equation}
Also, by \eqref{5.23-eq3},  $\tilde z\in L^2_{\dbF}(\O;H^1(0,T;L^2(G)))$ is the weak solution to the following random hyperbolic equation:
\begin{equation}\label{csz7.1}
\left\{\begin{array}{ll}
\ds \cA\tilde  z=\tilde  r_{1,t}+\tilde  r_2+\l ye^{2\l \phi}+\tilde  r &\mbox{ in } Q, \\[3mm]
\ns\ds \tilde  z=0   &\mbox{ on } \Si,\\[3mm]
\ns\ds  \tilde z(0)=\tilde z(T)=0  &\mbox{ in
}G.
\end{array}
\right.
\end{equation}
This implies  that
\begin{equation}\label{5.26-eq3}
\tilde z \in
L^2_\dbF(\Omega;C([0,T];H_0^1(\O)))\cap
L^2_\dbF(\Omega;C^1([0,T];L^2(\O))).
\end{equation}
The proof of (\ref{5.26-eq3}) is given in the
Appendix B.
For any constant $K>1$, put
 $$
 \tilde p\=K\tilde r.
 $$
By \eqref{5.23-eq4}-\eqref{5.24-eq7},  it is
easy to see that $\tilde p$ is the solution to
the following system:
\begin{equation}\label{csz7}
\left\{\begin{array}{ll} \ds
\cA \tilde p+\tilde {z}e^{-2\l\phi}=0    &\mbox{ in } Q,\\[2mm]
\ns\ds \tilde p = 0   &\mbox{ on } \Si,\\[2mm]
\ns\ds \tilde p(0)=\tilde p(T) =0  &\mbox{ in } G,\\[2mm]
\ns\ds \tilde p_t+
\varrho \frac{\tilde {r}_1}{\l^2}e^{-2\l \phi}=0    &\mbox{ in } Q,\\[2mm]
\ns\ds \tilde p-\varrho \frac{\tilde
{r}_2}{\l^4}e^{-2\l \phi}=0 &\mbox{ in } Q.
\end{array}
\right.
\end{equation}
Noting that $(\tilde r_1,\tilde
r_2)\in(L^2_{\dbF}(\O;H^1(0,T;L^2(G))))^2$,
similar to the proof of \eqref{5.26-eq3}, we can
also deduce that
$$
\tilde p\in
L^2_\dbF(\Omega;C([0,T];H_0^1(G)))\cap
L^2_\dbF(\Omega;C^1([0,T];L^2(G))).
$$

\smallskip

\noindent {\bf Step 2.}  Applying Theorem
\ref{5.24-lm1}  to $\tilde p$ in (\ref{csz7}), we obtain that
\begin{equation}\label{csz9}
\begin{array}{ll}
\ds \l\mE\int_Q(\l^2\tilde p^2+\tilde
p_t^2+|\nabla\tilde p|^2)e^{2\l \phi}dxdt
\\
\ns \ds\leq C\[\mE\int_Q\tilde {z}^2e^{-2\l
\phi}dxdt+\l^2\mE\int_0^T\int_{G_0}(\l^2
\tilde p^2+\tilde p_t^2)e^{2\l \phi}dxdt\]\\
\ns \ds \leq C\[\mE\int_Q\tilde {z}^2e^{-2\l
\phi}dxdt+\mE\int_0^T\int_{G_0}\(\frac{\tilde
{r}_1^2}{\l^2}+\frac{\tilde
{r}_2^2}{\l^4}\)e^{-2\l \phi}dxdt\].
\end{array}
\end{equation}
Here and hereafter,  $C$ denotes a constant,
independent of $K$ and $\l$.  Moreover, by (\ref{csz7})
again,   $\tilde p_t$ satisfies
\begin{equation}\label{csz10}
\left\{\begin{array}{ll} \ds
\cA \tilde p_t+(\tilde {z}e^{-2\l \phi})_t=0  &\mbox{ in  }\ Q,\\[2mm]
\ns\ds
\tilde p_t=0 &\mbox{ on }\ \Si,\\[2mm]
\ns\ds
\tilde p_{tt}+\frac{\varrho}{\l}\(\frac{\tilde r_{1,t}}{\l}-2\phi_t\tilde r_1\)e^{-2\l \phi}=0  &\mbox{ in }\  Q,\\[2mm]
\ns\ds \tilde p_t-\frac{\varrho
}{\l^2}\(\frac{\tilde
r_{2,t}}{\l^2}-\frac{2}{\l}\phi_t\tilde
r_2\)e^{-2\l \phi}=0  &\mbox{ in }\  Q.
\end{array}
\right.
\end{equation}
Applying Theorem \ref{5.24-lm1} to $\tilde p_t$,
 by (\ref{csz10}), we obtain that
\begin{equation}\label{5.24-eq10}
\begin{array}{ll}
\ds
\l\mE\int_Q\(\l^2\tilde p_t^2+\tilde p_{tt}^2+|\nabla\tilde p_t|^2\)e^{2\l \phi}dxdt \\
\ns \ds \leq \ds C\[\mE|e^{\l \phi}(e^{-2\l
\phi}\tilde {z})_t|^2_{L^2(Q)}
+\l^2\mE\int_0^T\int_{G_0}(\l^2\tilde p_t^2+\tilde p_{tt}^2)e^{2\l \phi}dxdt\]\\
\ns \ds \le \ds C\[\mE\int_Q(\tilde
z_t^2+\l^2\tilde
z^2)e^{-2\l \phi}dxdt +\mE\int_0^T\int_{G_0}\(\frac{\tilde
r_{1,t}^2}{\l^2}+\frac{\tilde
r_{2,t}^2}{\l^4}+\tilde r_1^2+\frac{\tilde
r_2^2}{\l^2}\)e^{-2\l \phi}dxdt\].
\end{array}
\end{equation}

\noindent {\bf Step 3.} By (\ref{csz7}), we have
that
$$
\ds -\mE\int_Q(\tilde r_{1,t}+\tilde
{r}_2)\tilde p dxdt =\mE\int_Q(\tilde r_1\tilde
p_t-\tilde r_2\tilde p)dxdt =-\mE\int_Q \varrho
\(\frac{\tilde r_1^2}{\l^2}+\frac{\tilde
r_2^2}{\l^4}\)e^{-2\l \phi}dxdt.
$$
This implies that
\begin{eqnarray*}
\begin{array}{ll}
\ds 0&=\mE(\cA\tilde {z}-\tilde r_{1,t}-\tilde
{r}_2-\l
y e^{2\l \phi}-\tilde  r,\tilde p)_{L^2(Q)}\\
\ns &\ds =-\mE\int_Q\tilde {z}^2e^{-2\l
\phi}dxdt-\mE\int_Q \varrho \(\frac{\tilde
r_1^2}{\l^2}+\frac{\tilde
r_2^2}{\l^4}\)e^{-2\l \phi}dxdt\\
\ns&\ds\q -\mE\l \int_Q y \tilde p e^{2\l
\phi}dxdt-K\mE\int_Q \tilde  r^2dxdt.
\end{array}
\end{eqnarray*}
Hence,
\begin{equation}\label{cspb7}
\mE\int_Q\tilde {z}^2e^{-2\l \phi}dxdt+\mE\int_Q
\varrho \(\frac{\tilde r_1^2}{\l^2}+\frac{\tilde
r_2^2}{\l^4}\)e^{-2\l \phi}dxdt +K\mE\int_Q
\tilde r^2dxdt  =-\l\mE\int_Q y \tilde p e^{2\l
\phi}dxdt.
\end{equation}
Combining  (\ref{csz9}) and (\ref{cspb7}), we
arrive at
\begin{equation}\label{5.23-eq01}
\mE\int_Q\tilde {z}^2e^{-2\l \phi}dxdt+\mE\int_Q
\varrho \(\frac{\tilde r_1^2}{\l^2}+\frac{\tilde
r_2^2}{\l^4}\)e^{-2\l \phi}dxdt  +K\mE\int_Q
\tilde r^2dxdt\leq \frac{C}{\l} \mE\int_Q y^2
e^{2\l \phi}dxdt.
\end{equation}

\noindent {\bf Step 4.} Using \eqref{csz7.1} and
\eqref{csz10} again, and noting $\tilde
p_{tt}(0)=\tilde p_{tt}(T)=0$ in $G$, we find
that
\begin{equation}\label{5.23-eq02}
\begin{array}{ll}
0\3n&\ds=\mE(\cA\tilde {z}-\tilde r_{1,t}-\tilde
{r}_2-\l
y e^{2\l \phi}-\tilde r,\tilde p_{tt})_{L^2(Q)}\\
\ns &\ds =-\mE\int_Q\tilde z (e^{-2\l
\phi}\tilde z)_{tt}dxdt-\mE\int_Q(\tilde
r_{1,t}+\tilde  r_2)\tilde p_{tt}dxdt\\
\ns&\ds\q-\l\mE\int_Q y\tilde p_{tt}e^{2\l
\phi}dxdt-\mE\int_Q\tilde r\tilde p_{tt}dxdt.
\end{array}
\end{equation}
Notice that
\begin{equation}\label{csj2}\begin{array}{ll}
\ds-\mE\int_Q\tilde z (e^{-2\l \phi}\tilde
z)_{tt}dxdt\\
\ns\ds =\mE\int_Q\[\tilde z_t^2e^{-2\l
\phi}dxdt-\frac{\tilde {z}^2}{2}
(e^{-2\l \phi})_{tt}\]dxdt\\
\ns  \ds=\mE\int_Q(\tilde z_t^2+\l
\phi_{tt}\tilde z^2-2\l^2\phi_t^2\tilde
z^2)e^{-2\l \phi}dxdt.
\end{array}
\end{equation}
Further, in view of the third and fourth
equalities in (\ref{csz10}), it follows that
\begin{equation}\label{csj3}
\begin{array}{ll}
\ds-\mE\int_Q(\tilde
r_{1,t}+\tilde  r_2)\tilde p_{tt}dxdt=-\mE\int_Q(\tilde r_{1,t}\tilde p_{tt}-\tilde  r_{2,t}\tilde p_t)dxdt\\
\ns \ds=\mE\int_Q \tilde
r_{1,t}\frac{\varrho}{\l}\(\frac{\tilde
r_{1,t}}{\l}-2\phi_t\tilde r_1\)e^{-2\l
\phi}dxdt +\mE\int_Q\tilde r_{2,t}\frac{\varrho
}{\l^2}\(\frac{\tilde
r_{2,t}}{\l^2}-\frac{2}{\l}\phi_t\tilde
r_2\)e^{-2\l \phi}dxdt\\
\ns \ds=\mE\int_Q\varrho \(\frac{\tilde
r_{1,t}^2}{\l^2}+\frac{\tilde
r_{2,t}^2}{\l^4}-\frac{2}{\l}\phi_t\tilde
r_1\tilde r_{1,t}-\frac{2}{\l^3}\phi_t\tilde
r_2\tilde r_{2,t}\)e^{-2\l \phi}dxdt.
\end{array}
\end{equation}
Further, by $\tilde p=K\tilde r$ and integration
by parts,
\begin{equation}\label{csj5}
\ds -\mE\int_Q\tilde r\tilde
p_{tt}dxdt=K\mE\int_Q\tilde r_t^2dxdt.
\end{equation}
Therefore, by (\ref{5.23-eq02})-(\ref{csj5}), we
get that
\begin{equation}\label{5.23-eq03}
\begin{array}{ll}
\ds \mE\int_Q\varrho \(\frac{\tilde
r_{1,t}^2}{\l^2}+\frac{\tilde
r_{2,t}^2}{\l^4}-\frac{2}{\l}\phi_t\tilde
r_1\tilde r_{1,t}-\frac{2}{\l^3}\phi_t\tilde
r_2\tilde r_{2,t}\)e^{-2\l
\phi}dxdt+K\mE\int_Q\tilde
r_t^2dxdt\\
\ns\ds +\mE\int_Q(\tilde z_t^2+\l
\phi_{tt}\tilde z^2-2\l^2\phi_t^2\tilde
z^2)e^{-2\l \phi}dxdt =\l\mE\int_Q y\tilde
p_{tt}e^{2\l \phi}dxdt.
\end{array}
\end{equation}

   Now, by (\ref{5.23-eq03}) and
(\ref{5.23-eq01}), using the Cauchy-Schwarz
inequality and noting (\ref{5.24-eq10}), we
obtain that
\begin{equation}\label{5.23-eq04}
\mE\int_Q(\tilde {z}_t^2\!+\!\l^2\tilde
{z}^2)e^{-2\l \phi}dxdt+\mE\int_Q \varrho
\(\frac{\tilde
{r}_{1,t}^2}{\l^2}\!+\!\frac{\tilde
{r}_{2,t}^2}{\l^4}\!+\!\tilde
r_1^2\!+\!\frac{\tilde r_2^2}{\l^2}\)e^{-2\l
\phi}dxdt  \!\leq\! C\l\mE\int_Q  y^2 e^{2\l
\phi}dxdt.
\end{equation}

\noindent {\bf Step 5.} By (\ref{csz7}),  we
find that
\begin{equation}\label{5.23-eq05}
\begin{array}{ll}
&\ds \mE(\tilde r_{1,t}+\tilde  r_2+\l ye^{2\l
\phi}+\tilde r,\tilde {z}e^{-2\l
\phi})_{L^2(Q)}=\mE(\cA\tilde  z, \tilde
{z}e^{-2\l \phi})_{L^2(Q)}\\
\ns&\ds =-\mE\int_Q\tilde z_t(\tilde {z}e^{-2\l
\phi})_tdxdt+\sum_{j,k=1}^n\mE\int_Qb^{jk}\tilde
z_{x_j}(\tilde {z}e^{-2\l \phi})_{x_k}dxdt\\
\ns&\ds =-\mE\int_Q(\tilde z_t^2+\l
\phi_{tt}\tilde
z^2-2\l^2\phi_t^2\tilde  z^2)e^{-2\l \phi}dxdt\\
\ns&\ds\q +\sum_{j,k=1}^n\mE\int_Qb^{jk}\tilde
z_{x_j}\tilde z_{x_k}e^{-2\l
\phi}dxdt-2\l\sum_{j,k=1}^n\mE\int_Qb^{jk}\tilde
z_{x_j}\tilde z\phi_{x_k}e^{-2\l\phi}dxdt.
\end{array}
\end{equation}
This yields that
\begin{equation}\label{csee15}
\begin{array}{ll}
\ds \mE\int_Q|\nabla\tilde
z|^2e^{-2\l\phi}dxdt\\
\ns \ds\le C\mE\int_Q\[|\tilde r_{1,t}+\tilde
r_2+\tilde r||\tilde {z}|e^{-2\l\phi}+\l|
y\tilde  z |+(\tilde
z_t^2+\l^2\tilde z^2)e^{-2\l \phi}\]dxdt\\
\ns \ds\le C\mE\int_Q\[y^2 e^{2\l
\phi}+\(\frac{\tilde
r_{1,t}^2}{\l^2}+\frac{\tilde
r_2^2}{\l^2}+\tilde r^2+\tilde z_t^2+\l^2\tilde
z^2\)e^{-2\l \phi}\]dxdt.
\end{array}
\end{equation}
By (\ref{5.23-eq01}), (\ref{5.23-eq04})
and (\ref{csee15}), we choose a constant $K$
in (\ref{5.23-eq01}) so that
\begin{equation}\label{zy01}
 K\ge Ce^{2\l\max\limits_{(t,x)\in
 Q}|\phi|}
 \end{equation}
(to absorb the term $C\ds\mE\int_Q\check
r^2e^{-2\l \phi}dxdt$ in (\ref{csee15})).  Then we deduce
that
\begin{equation}\label{5.23-eq06}
\begin{array}{ll}
\ds \mE\int_Q(|\nabla\tilde {z}|^2+\tilde
z_t^2+\l^2\tilde {z}^2)e^{-2\l \phi}dxdt
+\mE\int_Q \varrho \(\frac{\tilde
{r}_{1,t}^2}{\l^2}+\frac{\tilde
{r}_{2,t}^2}{\l^4}+\tilde r_1^2+\frac{\tilde
r_2^2}{\l^2}\)e^{-2\l \phi}dxdt
\\
\ns\ds\leq C\l\mE\int_Q y^2 e^{2\l \phi}dxdt.
\end{array}
\end{equation}

\noindent{\bf Step 6.} Recall that $(\tilde
{z},\tilde {r}_{1},\tilde {r}_{2},\tilde {r})$
depend on $K$. Now, we fix $\l$ and let $K$ tend
to infinity. By (\ref{5.23-eq01}) and
(\ref{5.23-eq06}),  we conclude that there
exists a subsequence of
$$
\begin{array}{ll}\ds
(\tilde {z},\tilde {r}_{1},\tilde {r}_{2},\tilde
{r}) \3n&\ds\in
(L^2_\dbF(\Omega;H^1_0(0,T;L^2(G)))\cap
L^2_\dbF(\Omega;L^2(0,T;H_0^1(G))))\\
\ns&\ds \q \times
(L^2_\dbF(\Omega;H^1(0,T;L^2(G))))^2\times
L^2_\dbF(0,T;L^2(G)),
\end{array}
$$
which converges weakly to some $(\bar
{z},\bar {r}_1,\bar {r}_2,0)$, with $\supp
\bar {r}_j\subset \cl{(0,T)\times G_0}$
($j=1,2$), since $\varrho(x)\equiv
\varrho^K(x)\to\infty$ for any $x\notin G_0$, as
$K\to\infty$. By (\ref{csz7}), we deduce that
$(\bar {z},\bar {r}_1,\bar {r}_2)$
satisfies
\begin{equation}\label{csxxz7}
\left\{\begin{array}{ll}
\ds \cA\bar  z=\bar  r_{1,t}+\bar  r_2+\l ye^{2\l \phi} &\mbox{ in }\ Q,\\[3mm]
\ns\ds \bar  z=0   &\mbox{ on }\ \pa Q.
\end{array}\right.
\end{equation}
Using (\ref{5.23-eq06}) again, we find that
\begin{equation}\label{5.24-eq8}
\mE|\bar {z}e^{-\l
\phi}|^2_{H_0^1(Q)}+\frac{1}{\l^2}\mE\int_0^T\int_{G_0}
(\bar {r}_{1,t}^2+\bar{r}_2^2)e^{-2\l
\phi}dxdt \leq C\l\mE\int_Q y^2 e^{2\l
\phi}dxdt.
\end{equation}
By (\ref{cspb6}), with $\eta$ replaced by
$\bar z$ above, one gets that
$$
\mE\(y,\bar r_{1,t}+\bar {r}_2+\l y e^{2\l
\phi}\)_{L^2(Q)}=\mE(b_1  y + f,\bar
{z})_{L^2(Q)}.
$$
Hence,  for any $\e>0$,
\begin{equation}\label{5.23-eq07}
\begin{array}{ll}
\ds \l \mE\int_Q  y^2 e^{2\l\phi}dxdt
\\[2mm]
\ns\ds =\mE(f, \bar {z})_{L^2(Q)} + \mE( b_1
y, \bar {z})_{L^2(Q)}-\mE( y,\bar
r_{1,t}+\bar  r_2)_{L^2((0,T)\times G_0)}\\[2mm]
\ns\ds \le C\Big\{\frac{1}{\e}\[ \mE|e^{\l\phi}
f|^2_{H^{-1}(Q)} + \mE|e^{\l
\phi}b_1 y|^2_{L^2(0,T;H^{-1}(G))}\\[2mm]
\ns\ds \q \q \q+\l^2\mE\int_0^T\int_{G_0}
y^2e^{2\l\phi}dxdt\]+\e\[\mE|e^{-\l \phi}\bar
z|^2_{H_0^1(Q)}+\mE|\bar
{z}e^{-\l\phi}|_{L^2(0,T;H^{1}_0(G))}^2
\\[2mm]
\ns\ds\q\q\q +\frac{1}{\l^2}\mE\int_0^T\int_{G_0}(\bar
r_{1,t}^2+\bar {r}_2^2)e^{-2\l
\phi}dxdt\]\Big\}.
\end{array}
\end{equation}
Finally, choosing $\e$ in (\ref{5.23-eq07})
sufficiently small and noting (\ref{5.24-eq8}),
 we arrive at the
desired estimate (\ref{5.24-eq9}). This
completes the proof of Theorem \ref{5.24-th1}.
\endpf

%%%%%%%%%%%%%%%%%%%%%%%%%%%%%%%%%%%%%%%%%%%%%%%%%%%%%

\section{An energy estimate for backward stochastic hyperbolic
equations}\label{sec-en-back}

%%%%%%%%%%%%%%%%%%%%%%%%%%%%%%%%%%%%%%%%%%%%%%%%%%%%%

In this section, we establish energy estimates respectively
for a random hyperbolic equation and a backward
stochastic hyperbolic equation, 
which will play important roles in the proof of
Theorem \ref{inobser}.

  First, set $\widehat T\in [0, T)$ and consider the following
random hyperbolic equation:
\begin{equation}\label{5.26-eq22}
\left\{
\begin{array}{ll}\ds
\vartheta_{tt}  -\sum_{i,j=1}^n
(b^{ij}\vartheta_{x_i})_{x_j}dt=  b_1\vartheta
&\mbox{ in }(\widehat T,T)\times G,\\[2mm]
\ns\ds \vartheta=0  &\mbox{ on }(\widehat T,T)\times \G,\\[2mm]
\vartheta(\widehat T)=\vartheta_0,\q
\vartheta_t(\widehat T)=\vartheta_1&\mbox{ in }G.
\end{array}
\right.
\end{equation}
It is easy to see that for  any $(\vartheta_0,
\vartheta_1)\in L^2(\Omega, \cF_{\widehat T},
\mathcal{P};  H^1_0(G))\times  L^2(\Omega,
\cF_{\widehat T}, \mathcal{P};  L^2(G))$,
\eqref{5.26-eq22} admits a unique solution
$$
\vartheta \in L^2_\mathbb{F}(\Omega; C([0, T];
H_0^1(G)))\bigcap L^2_\dbF(\O;C^1([0,
T];L^2(G))).
$$
Furthermore, we have the following energy estimate.
\begin{proposition}\label{prop2}
There is a constant $C>0$, depending only on
$T$, $G$ and $b^{ij}$ $(1\leq i,j\leq n)$, such
that for any solution $\vartheta$ to
\eqref{5.26-eq22} and for all $t,s$ satisfying $ \widehat
T\leq t \leq s \leq T$, it holds that
\begin{equation}\label{5.26-eq23}
\begin{array}{ll}
 \ds \mE\int_G  \big(|\vartheta_t(s,x)|^2 +
|\nabla\vartheta(s,x)|^2
\big) dx \\[3mm]
\ns \leq\ds Ce^{Cr_1^{\frac{1}{2-n/p}}(t-s)}
\mE\int_G \big(|\vartheta_t(t,x)|^2 + |\nabla
\vartheta(t,x)|^2 \big) dx.
\end{array}
\end{equation}
\end{proposition}

Next, let $\wt T\in (0,T]$. We consider the
following backward stochastic hyperbolic
equation:
\begin{equation}\label{5.25-eq17}
\left\{
\begin{array}{ll}\ds
d\a = \b dt + \eta dB(t) &\mbox{ in }(0,\wt T)\times G,\\
\ns\ds d\b-\sum_{i,j=1}^n
(b^{ij}\a_{x_i})_{x_j}dt=  b_1\a dt + b_2\eta dt
+ \zeta dB(t)
&\mbox{ in }(0,\wt T)\times G,\\[2mm]
\ns\ds \a=0,\ \b=0 &\mbox{ on }(0,\wt T)\times \G,\\[2mm]
\a(\wt T)=\a_0,\q \b(\wt T)=\b_0&\mbox{ in }G.
\end{array}
\right.
\end{equation}
Set
$$
\dbH_{\wt T}\!=\! L^2_\dbF(\O;C([0,\wt
T];H^1_0(G)))\!\times\! L^2_\dbF(\O;C([0,\wt
T];L^2(G))) \!\times\! L^2_\dbF(0,\wt
T;H^1_0(G))\!\times\! L^2_\dbF(0,\wt T;L^2(G)).
$$
We shall use  the following notion of solution for the
system (\ref{5.25-eq17}).

\begin{definition}\label{def sol-back}
$(\a,\b,\eta,\zeta)\in \dbH_{\wt T}$ is called a
solution to the system \eqref{5.25-eq17}, if
 \\[3mm]$(1)$ $\a(\wt T) = \a_0$ and
$\b(\wt T) = \b_0$ in $G$, $\mathcal{P}$-a.s.  \\[3mm]
$(2)$ For any $t \in (0,\wt T)$ and $\f \in
 C_0^1(G)$,  it holds that
\begin{equation} \label{solution to back1}
\begin{array}{ll}\ds
\q   \a(\wt T) -\a(t) =\int_t^{\wt T} \b(s)ds +
\int_t^{\wt T} \eta(s)dB(s) \q\mbox{ in }G,\q
\mathcal{P}\mbox{-a.s.}
\end{array}
\end{equation}
and
\begin{equation} \label{solution to back2}
\begin{array}{ll}\ds
  \langle \b(\wt T),\f\rangle_{L^2(G)} -
\langle \b(t),\f\rangle_{L^2(G)}
\\ \ns\ds= \int_t^{\wt T} \int_G \[
-\sum_{i,j=1}^n b^{ij}(x)\f_{x_i}(x)\a_{x_j}(s,
x) + b_1(s,x)\a(s,x) \f(x)\]dxds \\
\ns\ds \q   + \int_t^{\wt T} \int_G
\Big[b_2(s,x) \eta(s,x) \f(x) dxds+\zeta(s,
x)\f(x)dxdB(s)\],\q \mathcal{P}\mbox{-a.s.}
\end{array}
\end{equation}
\end{definition}

It is easy to show the following well-posedness
result for \eqref{5.25-eq17} (and hence we omit the proof).
\begin{lemma}\label{lm1} For any  $(\a_0, \b_0)\in  L^2(\Omega, \cF_T, \mathcal{P};  H^1_0(G))\times  L^2(\Omega, \cF_T, \mathcal{P};  L^2(G))$,
there is a unique solution
$(\a,\b,\eta,\zeta)\in \dbH_T$  to the system
\eqref{5.25-eq17}.
\end{lemma}

Furthermore, we have
the following energy estimate.
\begin{proposition}\label{prop1}
There is a constant $C>0$, depending only on
$T$, $G$ and $b^{ij}$ $(1\leq i,j\leq n)$, such
that for any solution $(\a,\b,\eta,\zeta)$ to
\eqref{5.25-eq17}, and for all $s,t$ satisfying  $0\leq s \leq t
\leq \wt T$, it holds that
\begin{eqnarray}\label{5.26-eq1}
\begin{array}{ll}
&\ds \mE\int_G  \big(|\b(s,x)|^2 + |\nabla
\a(s,x)|^2
\big) dx \\[3mm]
\ns& \leq\ds
Ce^{C\big(r_1^{\frac{1}{2-n/p}}+r_2^2\big)\wt
T} \mE\int_G \big(|\b(t,x)|^2 + |\nabla
\a(t,x)|^2 \big) dx,
\end{array}
\end{eqnarray}
and
\begin{eqnarray}\label{1012}
\begin{array}{rl}
&\ds\mE\int^{\wt T}_0 \int_G \Big(|\zeta|^2+|\nabla \eta|^2\Big)dxdt\\[3mm]
&\ds\leq
Ce^{C\big(r_1^{\frac{1}{2-n/p}}+r_2^2\big)\wt
T} |(\a_0, \b_0)|^2_{L^2(\Omega,
\cF_{\widetilde{T}}, \mathcal{P};
H^1_0(G))\times  L^2(\Omega,
\cF_{\widetilde{T}}, \mathcal{P};  L^2(G))}.
\end{array}\end{eqnarray}
\end{proposition}

\noindent {\bf Proof of Proposition \ref{prop1}.
} Define a (modified) energy of  the system
\eqref{5.25-eq17} as follows:
\begin{equation}\label{5.25-eq18}
\wt\cE(t) = \frac{1}{2}\mE\int_G \[|\b(t,x)|^2 +
\sum_{i,j=1}^n b^{ij}\a_{x_i}(t,x)\a_{x_j}(t,x)
+ r_1^{\frac{2}{2-n/p}}|\a(t,x)|^2\]dx,\q t\in
[0,T].
\end{equation}
Then, by It\^o's formula, we get that
\begin{equation}\label{5.25-eq19}
\begin{array}{ll}\ds
\wt\cE(t) -\wt \cE(s) \\
\ns\ds =  \mE\int_s^t\int_G \big[b_1\a\b +
b_2\b\eta \big]dxd\tau +
\frac{1}{2}\mE\int_s^t\int_G\Big(\sum\limits_{i, j=1}^{n} b^{i j}\eta_{x_i}\eta_{x_j}+ \zeta^2\Big) dxd\tau \\
\ns\ds \q+
r_1^{\frac{2}{2-n/p}}\mE\int_s^t\int_G \a\b
dxd\tau +
\frac{1}{2}r_1^{\frac{2}{2-n/p}}\mE\int_s^t\int_G
|\eta|^2 dxd\tau.
\end{array}
\end{equation}
Set $p_1 = \frac{2p}{n-2}$ and $p_2 =
\frac{2p}{p-n}$. By $\frac{1}{p} +
\frac{1}{p_1} + \frac{1}{p_2} + \frac{1}{2} = 1$
and $\frac{1}{2(n/p)^{-1}} +
\frac{1}{2(1-n/p)^{-1}} + \frac{1}{2} = 1$, and using
H\"older's inequality and Sobolev's embedding
theorem, we obtain that
$$
\begin{array}{ll}
\ds \Big|\mE\int_{G}b_1(\tau,x)
\a(\tau,x)\b(\tau,x)dx\Big|
\\[2mm]
\ns\ds \leq\mE\int_{G}|b_1(\tau,x)|
|\a(\tau,x)|^{\frac{n}{p}}|\a(\tau,x)|^{1-\frac{n}{p}}|\b(\tau,x)|dx\\[2mm]
\ns\ds \leq r_1\mE
\(\big||\a(\tau,\cd)|^{\frac{n}{p}}\big|_{L^{p_1}(G)}
\big||\a(\tau,\cd)|^{1-\frac{n}{p}}\big|_{L^{p_2}(G)}\big|
\b(\tau,\cd) \big|_{L^2(G)} \)\\[2mm]
\ns \ds = r_1 \mE \( \big|
\a(\tau,\cd)\big|^{\frac{n}{p}}_{L^{\frac{2n}{n-2}}(G)}
 \big|
\a(\tau,\cd)\big|^{1-\frac{n}{p}}_{L^{2}(G)}\big|\b(\tau,\cd)
\big|_{L^2(G)} \)\\[2mm]
\ns \ds = r_1^{\frac{1}{2-n/p}}\mE \( \big|
\a(\tau,\cd)\big|^{\frac{n}{p}}_{L^{\frac{2n}{n-2}}(G)}
r_1^{\frac{1-n/p}{2-n/p}}\big|
\a(\tau,\cd)\big|^{1-\frac{n}{p}}_{L^{2}(G)}\big|\b(\tau,\cd)
\big|_{L^2(G)} \).
\end{array}
$$
Notice that
$$
\!\left\{ \3n\begin{array}{ll}\ds \big|
\a(\tau,\cd)\big|^{\frac{n}{p}}_{L^{\frac{n}{n-2}}(G)}
\leq C\[\int_G\( |\b(\tau,x)|^2 + \sum_{i,j=1}^n
b^{ij}\a_{x_i}(\tau,x)\a_{x_j}(\tau,x) +
r_1^{\frac{2}{2-n/p}} |\a(\tau,x)|^2\) dx
\]^{\frac{n}{2p}},\\
\ns\ds r_1^{\frac{1-n/p}{2-n/p}}\big|
\a(\tau,\cd)\big|^{1-\frac{n}{p}}_{L^{2}(G)}\!\!
\leq\! \[\!\int_G\!\(\! |\b(\tau,x)|^2\! +\!\!
\sum_{i,j=1}^n\!
b^{ij}\a_{x_i}(\tau,x)\a_{x_j}(\tau,x)\! +\!
r_1^{\frac{2}{2-n/p}} |\a(\tau,x)|^2\!\) dx
\]^{\frac{1}{2}-\frac{n}{2p}},\\
\ns\ds \big|\b(\tau,\cd) \big|_{L^2(G)} \leq
\[\int_G\( |\b(\tau,x)|^2 + \sum_{i,j=1}^n
b^{ij}\a_{x_i}(\tau,x)\a_{x_j}(\tau,x) +
r_1^{\frac{2}{2-n/p}} |\a(\tau,x)|^2\) dx
\]^{\frac{1}{2}}.
\end{array}
\right.
$$
We have
\begin{equation}\label{en eq2.2}
\Big|\mE\int_{G}b_1(\tau,x)
\a(\tau,x)\b(\tau,x)dx\Big| \leq C
r_1^{\frac{1}{2-n/p}}\wt\cE(\tau).
\end{equation}
By a similar argument, we can also obtain that
\begin{equation}\label{en eq2.3}
\begin{array}{ll}\ds
 r_1^{\frac{2}{2-n/p}}\Big|\mE\int_G
\a(\tau,x)\b(\tau,x)dx\Big| \\
\ns\ds \leq
\frac{1}{2}r_1^{\frac{1}{2-n/p}}\mE\int_G\[
r_1^{\frac{2}{2-n/p}}\a^2(\tau,x) +
\b^2(\tau,x)\]dx \leq
r_1^{\frac{1}{2-n/p}}\wt\cE(\tau).
\end{array}
\end{equation}
Further, for a sufficiently small $\epsilon>0$,
\begin{equation}\label{5.26-eq2}
\begin{array}{ll}\ds
\Big|\mE\int_s^t\int_G
b_2(\tau,x)\b(\tau,x)\eta(\tau,x)
dxd\tau\Big|\\
\ns\ds \leq C(\epsilon)r_2^2
\int_s^t\wt\cE(\tau)d\tau  + \epsilon\mE\int_s^t\int_G|\nabla\eta(\tau,x)|^2
dxd\tau.
\end{array}
\end{equation}

By \eqref{5.25-eq19}-\eqref{5.26-eq2}, we find
that
\begin{equation}\label{en eq3}
\begin{array}{ll}
&\ds\frac{1}{2}\mE\int^t_s \int_G
\big(|\zeta|^2+r_1^{\frac{2}{2-n/p}}|\eta|^2+\frac{1}{2}
\sum\limits_{i, j=1}^{n} b^{i
j}\eta_{x_i}\eta_{x_j}\big)dxd\tau
+  \wt\cE(s) \\[3mm]
&  \leq \ds \wt\cE(t) +
C\(r_1^{\frac{1}{2-n/p}}+r_2^2\)
\int_s^t \wt\cE(\tau)d\tau.
\end{array}
\end{equation}
This, together with Gronwall's inequality,
implies that
\begin{equation}\label{en eq4}
\wt\cE(s) \leq
e^{C\big(r_1^{\frac{1}{2-n/p}}+r_2^2\big)\wt
T}\wt\cE(t),
\end{equation}
which implies \eqref{5.26-eq1}  and
(\ref{1012}).
\endpf

%%%%%%%%%%%%%%%%%%%%%%%%%%%%%%%%%%%%%%%%%%%%%%%%%%%%%

\section{Proof of Theorem \ref{inobser}}\label{sec-pr in}

%%%%%%%%%%%%%%%%%%%%%%%%%%%%%%%%%%%%%%%%%%%%%%%%%%%%%

In this section, we shall prove Theorem \ref{inobser}.

\ms

\noindent {\bf Proof of Theorem \ref{inobser}. }
We borrow some ideas from \cite{DZZ}.  The whole
proof  is divided into four steps.

\ms

\noindent {\bf Step 1.  }  Note that the
solution $y$ to \eqref{system1} may not be zero
at $t = 0$ and $t=T$. To apply Theorem
\ref{5.24-th1}, we need to choose a suitable
cutoff function.
Set
\begin{equation}\label{5.25-eq1}
\left\{
\begin{array}{ll}\ds
T_j = \frac{T}{2}-\e_j T,\q T_j'
=\frac{T}{2}+\e_j T,\\
\ns\ds R_0=\min_{x\in\overline G}\sqrt{d(x)},\q
R_1=\max_{x\in\overline G}\sqrt{d(x)},
\end{array}
\right.
\end{equation}
where $j=0,1,2$ and $0<\e_0<\e_1<\frac{1}{2}$.
By  \eqref{gz4}, (\ref{fllz101}) and (\ref{t0}), for any $T>T_0$, we have
that
\begin{equation}\label{5.25-eq2}
\phi(0,x)=\phi(T,x)\leq
R_1^2-\frac{c_1T^2}{4}<0,\q \forall\  x\in G.
\end{equation}
Therefore, there exists an $\e_1 \in (0,
\frac{1}{2})$,  which is close to $\frac{1}{2}$,
such that
\begin{equation}\label{5.25-eq3}
\phi(t,x) \leq
\frac{R_1^2}{2}-\frac{c_1T^2}{8}<0,\q \forall\
(t,x)\in [(0,T_1)\cup(T_1',T_1)]\times G.
\end{equation}

On the other hand, it follows from \eqref{gz4}
that
$$
\phi\(\frac{T}{2},x\)=d(x)\geq R_0^2,\q \forall\
x\in G.
$$
Therefore, there is an $\e_0\in (0,
\frac{1}{2})$, which is close to $0$, such that
\begin{equation}\label{5.25-eq4}
\phi(t,x) \geq \frac{R_0^2}{2},\q \forall\
(t,x)\in (T_0,T_0')\times G.
\end{equation}
Furthermore,  choose a nonnegative function
$\xi\in C_0^\infty(0, T)$ such that
\begin{equation}\label{5.25-eq5}
\xi(t)=1 \q\mbox{ in }(T_1,T_1').
\end{equation}

\noindent {\bf Step 2.} In this step, we prove
that there is a $\l_1>0$, such that for any
$\l\geq
\l_1$,
\begin{equation}\label{5.25-eq6}
\l\mE\int_Q e^{2\lambda\phi} y^2 dxdt \leq
C\(\l^2\mE\int_0^T\int_{G_0}e^{2\lambda\phi} y^2 dxdt +
\mE|y|_{L^2(J\times G)}^2+\mE\int_Q e^{2\lambda\phi}
f^2dxdt\),
\end{equation}
where $J = (0, T_1)\cup (T_1', T)$.

To this aim, set $\tilde y=\xi y$.  Then $\tilde
y$  satisfies the following forward stochastic
hyperbolic equation:
\begin{eqnarray}\label{1013}
\left\{
\begin{array}{lll}\ds
\ds d\tilde{y}_{t} -
\sum_{i,j=1}^{n}(b^{ij}\tilde{y}_{x_i})_{x_j}dt
= \big( b_1 \tilde{y} + \tilde f \big)dt +
\big(b_2 \tilde{y} + \xi g\big)dB(t) &{\mbox {in
}} Q,
 \\
\ns\ds  \tilde{y} = 0  &\mbox{on } \Si, \\[3mm]
\ns\ds  \tilde{y}(0) =\tilde{y}(T)=0 &\mbox{in }
G,
\end{array}
\right.
\end{eqnarray}
with $\tilde f=\xi f+\xi_{tt}y+2\xi_t y_t$. By
Theorem \ref{5.24-th1}, for any $\l\geq\l_0$, we
have that
\begin{equation}\label{5.25-eq7}
\begin{array}{ll}\ds
\l\mE\int_Q e^{2\lambda\phi} \tilde y^2 dxdt\\ \ns \ds \leq
C\(\mE|e^{\lambda\phi} \tilde f|^2_{H^{-1}(Q)}+
|e^{\lambda\phi} b_1 \tilde
y|^2_{L^2_\dbF(0,T;H^{-1}(G))}+\l^2
\mE\int_0^T\int_{G_0} e^{2\lambda\phi} \tilde y^2
dxdt\).
\end{array}
\end{equation}

By the definition of $\tilde f$, we find that
\begin{eqnarray}\label{1015}
\begin{array}{rl}
&\mE|e^{\lambda\phi} \tilde f|^2_{H^{-1}(Q)}=
\mE|e^{\lambda\phi} (\xi f+\xi_{tt}y+2\xi_t y_t)|^2_{H^{-1}(Q)}\\[3mm]
&=\sup\limits_{|h|_{L^2(\Omega; H^1_0(Q))}=1}
\Big|\mE \langle  e^{\lambda\phi} (\xi f+\xi_{tt}y+2\xi_t y_t), h\rangle_{H^{-1}(Q), H^1_0(Q)}\Big|^2\\[3mm]
&\leq |e^{\lambda\phi} f|^2_{L^2_\mathbb{F} (0, T;
L^2(G))}+C\l^2
|e^{\lambda\phi} y|^2_{L^2_\mathbb{F}(J; L^2(G))}\\[3mm]
&\leq |e^{\lambda\phi} f|^2_{L^2_\mathbb{F} (0, T;
L^2(G))}+C\l^2 e^{\l(R_1^2-\frac{c_1T^2}{4})}
|y|^2_{L^2_\mathbb{F}(J; L^2(G))}.
\end{array}
\end{eqnarray}
Further, recalling the definition of $r_1$ and
noting the embedding $L^{2p/p+2}(G)
\hookrightarrow H^{-1}(G)$, we get that
\begin{equation}\label{5.25-eq9}
|e^{\lambda\phi}b_1 \tilde
y|_{L^2_\dbF(0,T;H^{-1}(G))}\leq  C|e^{\lambda\phi}b_1
\xi y|_{L^2_\dbF(0,T;L^{2p/p+2}(G))}\leq
Cr_1|e^{\lambda\phi}\xi y|_{L^2_\dbF(0,T;L^{2}(G))}.
\end{equation}
Further, by \eqref{5.25-eq3} and
\eqref{5.25-eq5},
\begin{equation}\label{5.25-eq10}
\begin{array}{ll}\ds
|e^{\lambda\phi}\xi
y|_{L^2_\dbF(0,T;L^{2}(G))}^2 \3n& \ds=
|e^{\lambda\phi} y|_{L^2_\dbF(0,T;L^{2}(G))}^2 -
\mE\int_Q e^{2\lambda\phi}(1-\xi^2)y^2 dxdt \\
\ns&\ds \geq |e^{\lambda\phi}
y|_{L^2_\dbF(0,T;L^{2}(G))}^2 -
Ce^{(R_1^2-c_1T^2/4)\l}|y|^2_{L^2_\mathbb{F}(J;
L^2(G))}.
\end{array}
\end{equation}
Therefore, by
\eqref{5.25-eq7}-\eqref{5.25-eq10},    there is
a constant $C_1 = C_1(T, G)$, independent of
$\l$ and $r_1$, such that
\begin{equation}\label{5.25-eq11}
\begin{array}{ll}\ds
|e^{\lambda\phi}
y|_{L^2_\dbF(0,T;L^{2}(G))}^2\3n&\ds\leq C_1\[
\frac{r_1^2}{\l}|e^{\lambda\phi}
y|_{L^2_\dbF(0,T;L^{2}(G))}^2 +
\l\mE\int_0^T\int_{G_0}e^{2\lambda\phi} y^2 dxdt \\
\ns&\ds\qq \ +
e^{(R_1^2-cT^2/4)\l}(1+\l)\mE|y|^2_{L^2(J\times
G)}+\frac{1}{\l}|e^{\lambda\phi}
f|^2_{L^2_{\mathbb{F}}(0, T; L^2(G))}\].
\end{array}
\end{equation}
Since $R^2_1-cT^2/4 < 0$, one may find a
sufficiently large $\l_1>0$,
 such that for any $\l>\l_1$, (\ref{5.25-eq6})
 holds.

\medskip

\noindent {\bf Step 3. }  We establish an energy
estimate for solutions to (\ref{system1}). Set
\begin{equation}\label{5.25-eq14}
\cE(t)\=\frac{1}{2}\big(\mE|y(t,\cd)|_{L^2(G)}^2
+ \mE|y_t(t,\cd)|_{H^{-1}(G)}^2 \big).
\end{equation}
Then by the classical energy estimate, for any
$S_0\in (T_0,\frac{T}{2})$ and $S_0'\in
(\frac{T}{2},T_0')$,
\begin{equation}\label{5.25-eq15}
\int_{S_0}^{S_0'}\cE(t)dt\leq C(1+r_1 +
r_2)\mE\int_{S_0}^{S_0'}\int_G y^2dxdt +
C\mE\int_0^T\int_G (f^2+g^2)dxdt.
\end{equation}

On the other hand,  we claim that there exists
a constant $C > 0$,  such that
\begin{equation}\label{5.25-eq16}
\cE(t)\leq
Ce^{C\big(r_1^{\frac{1}{2-n/p}}+r_2^2\big)}\big(\cE(s)+|(f,
g)|_{(L^2_\mathbb{F}(0, T;
L^2(G)))^2}^2\big),\qq \forall \  t,\,s\in
[0,T].
\end{equation}
In the following,  we only prove the case of
$t\geq s$.  The  other case can be also proved
by a similar technique and Proposition
\ref{prop2}.
 By It\^o's formula,  let $\widetilde T=t$ in
\eqref{5.25-eq17} and $T=t$ in (\ref{system1}).
Then it follows that
\begin{equation}\label{5.25-eq20}
\begin{array}{ll}\ds
\mE\langle y(t),\b_0 \rangle_{L^2(G)} +
\mE\langle
y_t(t),-\a_0\rangle_{H^{-1}(G),H_0^1(G)}\\[3mm]
\ns\ds = \mE\langle y(s),\b(s) \rangle_{L^2(G)}
+ \mE\langle
y_t(s),-\a(s)\rangle_{H^{-1}(G),H_0^1(G)}\\[3mm]
\ns\ds \quad-\mE\int^t_s\int_G  (\a f+\eta
g)dxdt.
\end{array}
\end{equation}
Denote by $\dbS$ the unit sphere of the space
$L^2_{\cF_t}(\O;H_0^1(G))\times
L^2_{\cF_t}(\O;L^2(G))$. By  \eqref{5.25-eq20},
 \eqref{5.26-eq1} and \eqref{1012}, we obtain that
$$
\begin{array}{ll}&\ds
\sqrt{2\cE(t)} \ds= \sup_{(\a_0,\b_0)\in
\dbS}\big|\mE\big(\langle y(t),\b_0
\rangle_{L^2(G)} + \langle
y_t(t),-\a_0\rangle_{H^{-1}(G),H_0^1(G)}\big)\big|\\
\ns&\ds  = \sup_{(\a_0,\b_0)\in
\dbS}\Big|\mE\[\langle y(s),\b(s)
\rangle_{L^2(G)} + \langle
y_t(s),-\a(s)\rangle_{H^{-1}(G),H_0^1(G)}-\int^t_s\int_G  (\a f+\eta g)dxdt\]\Big|\\
\ns&\ds \leq C\sqrt{\cE(s)}\sup_{(\a_0,\b_0)\in
\dbS}|(\a(s),\b(s))|_{L^2_{\cF_s}(\O;H_0^1(G))\times
L^2_{\cF_s}(\O;L^2(G))}\\
\ns&\ds\quad+\sup_{(\a_0,\b_0)\in
\dbS}|(\a, \eta)|_{L^2_\mathbb{F}(s, t; L^2(G))}|(f, g)|_{(L^2_\mathbb{F}(s, t; L^2(G)))^2}\\
\ns&\ds \leq
Ce^{C\big(r_1^{\frac{1}{2-n/p}}+r_2^2\big)}\big(\sqrt{\cE(s)}+|(f,
g)|_{(L^2_\mathbb{F}(s, t;
L^2(G)))^2}\big).\end{array}
$$
This implies our claim (\ref{5.25-eq16}).

\medskip

\noindent {\bf Step 4}. First, it follows from
\eqref{5.25-eq4}  that
\begin{equation}\label{5.25-eq13}
\mE\int_Q e^{2\lambda\phi} y^2 dxdt \geq
e^{R_0^2\l}\mE\int_{T_0}^{T_0'}\int_G y^2 dxdt.
\end{equation}
Also,  \eqref{5.25-eq16} implies that
\begin{equation}\label{5.25-eq21}
\mE|y|_{L^2(J\times G)}^2\leq
Ce^{C\big(r_1^{\frac{1}{2-n/p}}+r_2^2\big)}\big(\cE(0)+|(f,
g)|_{(L^2_\mathbb{F}(0, T; L^2(G)))^2}^2\big)
\end{equation}
and
\begin{equation}\label{5.25-eq22}
Ce^{-C\big(r_1^{\frac{1}{2-n/p}}+r_2^2\big)}\cE(0)
\leq \int^{S'_0}_{S_0} \cE(t)dt+|(f,
g)|_{(L^2_\mathbb{F}(0, T; L^2(G)))^2}^2.
\end{equation}
By \eqref{5.25-eq13}, \eqref{5.25-eq15} and
\eqref{5.25-eq22},  we get that
\begin{equation}\label{5.25-eq23}
\cE(0)\leq
Ce^{C\big(r_1^{\frac{1}{2-n/p}}+r_2^2\big)}
\big(e^{-R_0^2\l}\mE\int_Q e^{2\lambda\phi}
y^2dxdt+|(f, g)|_{(L^2_\mathbb{F}(0, T;
L^2(G)))^2}^2\big).
\end{equation}
Combining the above estimate with
(\ref{5.25-eq6}) and  (\ref{5.25-eq21}), we find
that
$$
\begin{array}{ll}\ds
\cE(0)\3n&\ds\leq
Ce^{C\big(r_1^{\frac{1}{2-n/p}}+r_2^2\big)}
\Big[e^{-R_0^2\l}\Big(\l\mE\int^T_0\int_{G_0} e^{2\lambda\phi}y^2dxdt+\frac{1}{\l}\mE|y|^2_{L^2(J\times G)}\Big)\\
\ns&\ds\qq\qq\qq\qq+ |(f, g)|_{(L^2_\mathbb{F}(0, T; L^2(G)))^2}^2\Big]\\
\ns&\ds\leq
Ce^{C\big(r_1^{\frac{1}{2-n/p}}+r_2^2\big)}
\Big[e^{-R_0^2\l}\Big(\l\mE\int^T_0\int_{G_0} e^{2\lambda\phi}y^2dxdt+\frac{1}{\l}\cE(0)\Big)\\
\ns&\ds\qq\qq\qq\qq+ |(f,
g)|_{(L^2_\mathbb{F}(0, T; L^2(G)))^2}^2\Big].
\end{array}
$$
Therefore, there exists a sufficiently large
constant $\l_3>0$, such that for any $\l>\l_3$,
 the
desired observability inequality \eqref{inobser
esti2.1} holds.
\endpf

%%%%%%%%%%%%%%%%%%%%%%%%%%%%%%%%%%%%%%%%%%%%%%%%%%%%%

\section{Appendix A: Proof of Proposition \ref{5.24-prop1}}

In this section, we give a proof of  Proposition
\ref{5.24-prop1}. To this aim, we need
the following known result.

\begin{lemma}\label{5.24-lm2} $\cite[Proposition\  3.5]{Fu-Yong-Zhang1}$
For any $h>0$,  $m=3, 4, \cdots$, and
$q_m^j,w_m^j\in \dbC$ $(j=0,1,\cds,m)$ satisfying
$q_m^0=q_m^m=0$,  we have that
\begin{equation}\label{5.23-eq2}
\begin{array}{ll}\ds
-\sum_{j=1}^{m-1}q_m^j\frac{w_m^{j+1}-2w_m^j+w_m^{j-1}}{h^2}\3n&\ds=\sum_{j=0}^{m-1}\frac{q_m^{j+1}-q_m^j}{h}\frac{w_m^{j+1}-w_m^j}{h}\\
\ns
&\ds=\sum_{j=1}^m\frac{q_m^j-q_m^{j-1}}{h}\frac{w_m^j-w_m^{j-1}}{h}.
\end{array}
\end{equation}
\end{lemma}

\noindent {\bf Proof of Proposition
\ref{5.24-prop1}.}  The whole proof is divided
into four steps.

\medskip

\noindent {\bf Step 1.}  Let
$\Big\{\{(z_m^{j,k},r_{1m}^{j,k}, r_{2m}^{j,k},
r_m^{j,k})\}_{j=0}^m\Big\}_{k=1}^\infty\subset\cA_{ad}$
be  a minimizing sequence of $J(\cd)$. Thanks to
the coercivity of the cost functional, since
$z_m^{j,k}$ solves an elliptic equation, it can
be shown that $\Big\{\{(z_m^{j,k},r_{1m}^{j,k},
r_{2m}^{j,k},
r_m^{j,k})\}_{j=0}^m\Big\}_{k=1}^\infty$ is
bounded in $\cA_{ad}$. Therefore, there exists a
subsequence of $\Big\{\{(z_m^{j,k},r_{1m}^{j,k},
r_{2m}^{j,k},
r_m^{j,k})\}_{j=0}^m\Big\}_{k=1}^\infty$
converging  weakly  to some $\{(\hat z_m^j,\hat
r_{1m}^j, \hat r_{2m}^j, \hat
r_m^j)\}_{j=0}^m\in \cA_{ad}$ in
$(L^2_{\cF_{jh}}(\O;H_0^1(G))\times
(L^2_{\cF_{jh}}(\O;L^2(G)))^3)^{m+1}$. Since the
functional $J$ is strictly convex, this element
is the unique solution to (\ref{5.23-eq1}). By
(\ref{csxg1}) and the definition of
 $\cA_{ad}$, it is obvious that $\hat z_m^0=\hat
z_m^m=p_m^0=p_m^m=0$ in $G$.

\medskip

\noindent {\bf Step 2.} Fix  $\d_{0m}^j\in
L^2_{\cF_{jh}}(\O;H^2(G)\cap H_0^1(G))$,
$\d_{1m}^j\in L^2_{\cF_{jh}}(\O;L^2(G))$, and
$\d_{2m}^j\in L^2_{\cF_{jh}}(\O;L^2(G))$
($j=0,1,2,\cdots,m$) with
$\d_{0m}^0=\d_{0m}^m=\d_{2m}^0=\d_{2m}^m\equiv
0$ and $\d_{1m}^0=\d_{1m}^1$. For
$(\mu_0,\mu_1,\mu_2)\in\dbR^3$, put
$$
\left\{
\begin{array}{ll}
\ds r_m^j\=\mE\(\frac{\hat z_m^{j+1}-2\hat
z_m^j+\hat z_m^{j-1}}{
h^2}\ \!\Big|\ \!\cF_{jh}\)+\mE\(\frac{\d_{0m}^{j+1}-2\d_{0m}^j+\d_{0m}^{j-1}}{
h^2}\ \!\Big|\ \!\cF_{jh}\)\mu_0\\
\ds\qq\;\;
-\sum_{j_1,j_2=1}^n\pa_{x_{j_2}}\(b^{j_1j_2}\pa_{x_{j_1}}(\hat
z_m^j+\mu_0\d_{0m}^j)\)-\mE\(\frac{\hat
r_{1m}^{j+1}-\hat r_{1m}^j}{
h}\ \!\Big|\ \!\cF_{jh}\)\\
\ns
\qq\;\;\ds-\mE\(\frac{\d_{1m}^{j+1}-\d_{1m}^j}{
h}\ \!\Big|\ \!\cF_{jh}\)\mu_1-\hat
r_{2m}^j-\mu_2\d_{2m}^j
-\l y_m^je^{2\l \phi_m^j},\qq1\le j\le m-1;\\[3mm]
\ns r_m^0=r_m^m=0.
\end{array}
\right.
$$Then  $\{(\hat
z_m^j+\mu_0\d^j_{0m},\hat
r_{1m}^j+\mu_1\d^j_{1m},\hat
r_{2m}^j+\mu_2\d^j_{2m},r^j_m)\}_{j=0}^m\in\cA_{ad}$.
Define a function $g(\cd,\cd,\cd)$ in $\dbR^3$
by
$$
g(\mu_0,\mu_1,\mu_2)=J\big(\{(\hat
z_m^j+\mu_0\d_{0m}^j,\hat
r_{1m}^j+\mu_1\d_{1m}^j,\hat
r_{2m}^j+\mu_2\d_{2m}^j, r_m^j)\}_{j=0}^m\big).
$$
Since $g$ has a minimum at $(0,0,0)$,
we get that  $\nabla g(0,0,0)=0$.

\medskip

By $\ds\frac{\pa g(0,0,0)}{\pa\mu_1}=\ds\frac{\pa
g(0,0,0)}{\pa\mu_2}=0$,  noting that $\{(\hat
z_m^j,\hat r_{1m}^j, \hat r_{2m}^j, \hat
r_m^j)\}_{j=0}^m$ satisfy the first equation of
(\ref{5.23-eq3}), we find that
$$
\begin{array}{ll}
\ds-K\mE\sum_{j=1}^{m-1}\int_G \hat
r_m^j\mE\(\frac{\d_{1m}^{j+1}-\d_{1m}^j}{
h}\ \!\Big|\ \!\cF_{jh}\)dx+\mE\sum_{j=1}^m\int_G
\varrho \frac{\hat r_{1m}^j\d_{1m}^j}{\l^2}e^{-2\l \phi_m^j}dx\\
\ns\ds = K\mE\sum_{j=1}^{m}\int_G \frac{\hat
r_m^j-\hat r_m^{j-1}}{
h}\d_{1m}^jdx+\mE\sum_{j=1}^m\int_G \varrho
\frac{\hat r_{1m}^j\d_{1m}^j}{\l^2}e^{-2\l
\phi_m^j}dx \\
\ns\ds = K\mE\sum_{j=1}^{m}\int_G \[\frac{\hat
r_m^j-\hat r_m^{j-1}}{ h}+\varrho \frac{\hat
r_{1m}^j}{\l^2}e^{-2\l \phi_m^j}\]\d_{1m}^jdx
=0,
\end{array}
$$

and
$$
\begin{array}{ll}\ds
-K\mE\sum_{j=1}^{m-1}\int_G \hat
 r_m^j\d_{2m}^jdx+\mE\sum_{j=1}^{m-1}\int_G
\varrho \frac{\hat r_{2m}^j\d_{2m}^j}
{\l^4}e^{-2\l \phi_m^j}dx\\
\ns\ds=-\mE\sum_{j=1}^{m-1}\int_G \(K
\hat r_m^j-\varrho \frac{\hat r_{2m}^j}
{\l^4}e^{-2\l \phi_m^j}\)\d_{2m}^jdx=0,
\end{array}
$$
which, combined with (\ref{csxg1}),  gives
(\ref{5.23-eq4}).

\medskip

On the other hand, by $\ds\frac{\pa
g(0,0,0)}{\pa\mu_0}=0$, we have that
\begin{eqnarray}\label{csg3}
\begin{array}{rl}
&\ds\mE\!\sum_{j=1}^{m-1}\!\int_G\! \Big\{K\hat
r_m^j\[\mE\(\frac{\d_{0m}^{j+1}\!\!-2\d_{0m}^j\!+\d_{0m}^{j-1}}{
h^2}\ \!\Big|\ \!\cF_{jh}\)-\!\!\sum_{j_1,j_2=1}^n\!\!\!\pa_{x_{j_2}}\!(b^{j_1j_2}\pa_{x_{j_1}}\d_{0m}^j)\]\\[3mm]
&\quad\quad\quad\quad\quad+\!\hat
z_m^j\d_{0m}^je^{-2\l \phi_m^j}\Big\}dx=0,
\end{array}
\end{eqnarray}
which, combined with
$p_m^0=p_m^m=\d_{0m}^0=\d_{0m}^m=0$ in $G$,
implies that (\ref{cs112}) holds. By means of
the regularity theory for elliptic equations of
second order, one finds that $\hat
z_m^j,\; p_m^j\in
L^2_{\cF_{jh}}(\O;H^2(G)\cap H_0^1(G))$, $1\leq
j\leq m-1$.

\medskip

\noindent{\bf Step 3.} Recalling that $\{(\hat
z_m^j,\hat r_{1m}^j, \hat r_{2m}^j, \hat
r_m^j)\}_{j=0}^m$ satisfy (\ref{5.23-eq3}), and
noting (\ref{5.23-eq4})-(\ref{cs112}) and
$p_m^j=K\hat r_m^j$, one gets
\begin{equation}\label{5.23-eq5}
\begin{array}{ll}
\ds 0=\mE\sum_{j=1}^{m-1}\int_G
\[\mE\(\frac{\hat z_m^{j+1}-2\hat z_m^j+\hat
z_m^{j-1}}{
h^2}\ \!\Big|\ \!\cF_{jh}\)-\sum_{j_1,j_2=1}^n\pa_{x_{j_2}}(b^{j_1j_2}\pa_{x_{j_1}}\hat z_m^j)\\
\ns\ds\qq\qq\qq-\mE\(\frac{\hat
r_{1m}^{j+1}-\hat
r_{1m}^j}{h}\ \!\Big|\ \!\cF_{jh}\)-\hat r_{2m}^j
-\l y_m^je^{2\l \phi_m^j}-\hat r_m^j\]p_m^jdx\\
\ns\ds \q=\mE\sum_{j=1}^{m-1}\int_G
\[\mE\(\frac{p_m^{j+1}-2p_m^j+p_m^{j-1}}{h^2}\ \!\Big|\ \!\cF_{jh}\)-\sum_{j_1,j_2=1}^n\pa_{x_{j_2}}(b^{j_1j_2}\pa_{x_{j_1}}p_m^j)\]\hat z_m^jdx\\
\ns\ds\qq +\mE\sum_{j=1}^m\int_G
\frac{p_m^j-p_m^{j-1}}{h}\hat
r_{1m}^jdx-\mE\sum_{j=1}^{m-1}\int_G\(\hat
r_{2m}^j +\l y_m^je^{2\l \phi_m^j}+\hat
r_m^j\)p_m^jdx
\\
\ns\ds\q =-\mE\sum_{j=1}^{m-1}\[\int_G|\hat
z_m^j|^2 e^{-2\l \phi_m^j}dx +\int_G \varrho
\(\frac{|\hat r_{1m}^j|^2}{\l^2}+\frac{|\hat
r_{2m}^j|^2}{\l^4}\)e^{-2\l \phi_m^j}dx \\
\ns\ds\qq +K\int_G|\hat r_m^j|^2 dx\]-\mE\int_G
\varrho \frac{|\hat r_{1m}^m|^2}{\l^2}
e^{-2\lambda \phi^m_m}dx-\l
\mE\sum_{j=1}^{m-1}\int_G y_m^je^{2\l
\phi_m^j}p_m^jdx.
\end{array}
\end{equation}
By (\ref{5.23-eq4}) and (\ref{5.23-eq5}),   there is a constant $C=C(K,
\lambda)>0$, such that
$$
\begin{array}{ll}\ds
\mE\sum_{j=1}^{m-1}\[\int_G |\hat z_m^j|^2
e^{-2\l \phi_m^j}dx
 +\int_G \varrho \(\frac{|\hat
r_{1m}^j|^2}{\l^2}+\frac{|\hat
r_{2m}^j|^2}{\l^4}\)e^{-2\l \phi_m^j}dx\\
\ns\ds\qq +K\int_G|\hat r_m^j|^2 dx\]+\mE\int_G
\varrho \frac{|\hat
r_{1m}^m|^2}{\l^2}e^{-2\lambda \phi^m_m}dx \leq
C\mE\sum_{j=1}^{m-1}\int_G|y_m^j|^2e^{2\l
\phi_m^j}dx.
\end{array}
$$
This implies (\ref{5.23-eq6}).

\medskip

\noindent {\bf Step 4.}  Noting that (\ref{cs112}) holds  and
$p_m^0=\hat z_m^0=p_m^m=\hat z_m^m=0$, we obtain that
\begin{equation}\label{csOK22}
\3n\3n\begin{array}{ll} \ds
\frac{\mE(p_m^3\ \!|\ \!\cF_{h})
-4\mE(p_m^2\ \!|\ \!\cF_{h})+5p_m^1}{
h^4} -\sum_{j_1,j_2=1}^n\pa_{x_{j_2}}\[b^{j_1j_2}\pa_{x_{j_1}}\mE\(\frac{p_m^2-2p_m^1+p^0_m}{h^2}\ \!\Big|\ \!\cF_{h}\)\]\\[3mm]
\ns \ds\qq\q+\frac{\mE(\hat z_m^2\ \!|\ \!\cF_h)e^{-2\l
\phi_m^2}-2\hat z_m^1e^{-2\l \phi_m^1}+\hat
z^0_me^{-2\l \phi^0_m}}{h^2}=0\qq \mbox{ in }G,\\ \ns
\ds  \frac{4p_m^{m-1}\!
+\!\mE(p_m^{m-1}\ \!|\ \!\cF_{(m-2)h})\!-\!4p_m^{m-2}\!+\!p_m^{m-3}}{
h^4}\\ \ns
\qq\q\ds-\!\sum_{j_1,j_2=1}^n\!\!\pa_{x_{j_2}}\!\!\[b^{j_1j_2}\pa_{x_{j_1}}\!\mE\(\frac{p_m^m\!-\!2
p_m^{m-1}\! +\! p_m^{m-2} }{
h^2}\ \!\Big|\ \!\cF_{(m-1)h}\)\]\\[3mm]
\ns \ds\qq\q+\frac{\hat z^m_me^{-2\l
\phi_m^m}-2\hat z_m^{m-1} e^{-2\l
\phi_m^{m-1}}+\hat z_m^{m-2} e^{-2\l
\phi_m^{m-2}}}{h^2}=0\qq \mbox{ in }\ G,
\end{array}
\end{equation}
and for $j=2,\cdots,m-2$,
\begin{equation}\label{cs22}
\begin{array}{ll}
\ds \frac{\mE(p_m^{j+2}\ \!|\ \!\cF_{jh})
-4\mE(p_m^{j+1}\ \!|\ \!\cF_{jh})+5p_m^j+\mE(p_m^j\ \!|\ \!\cF_{(j-1)h})-4p_m^{j-1}+p_m^{j-2}}{
h^4} \\
\ns\ds
-\sum_{j_1,j_2=1}^n\pa_{x_{j_2}}\[b^{j_1j_2}\pa_{x_{j_1}}\mE\(\frac{p_m^{j+1}-2p_m^j+p_m^{j-1}}{
h^2}\ \!\Big|\ \!\cF_{jh}\)\]\\
\ns \ds +\frac{\mE(\hat z_m^{j+1}\ \!|\ \!\cF_{jh})
e^{-2\l \phi_m^{j+1}}-2\hat z_m^j e^{-2\l
\phi_m^j}+\hat z_m^{j-1} e^{-2\l \phi_m^{j-1}}}{
h^2}=0\q \mbox{ in }G.
\end{array}
\end{equation}
By \eqref{5.23-eq3},
\begin{equation}\label{cs33}
\begin{array}{ll}
0\3n &\ds=\mE\sum_{j=1}^{m-1}\int_G
\[\mE\(\frac{\hat
z_m^{j+1} -2\hat z_m^j +\hat
z_m^{j-1}}{h^2}\ \!\Big|\ \!\cF_{jh}\)
 - \sum_{j_1,j_2=1}^n
\pa_{x_{j_2}}\big(b^{j_1j_2}\pa_{x_{j_1}}\hat
z_m^j\big)\\
\ns&\ds \q - \mE\(\frac{\hat r_{1m}^{j+1} - \hat
r_{1m}^j}{h}\ \!\Big|\ \!\cF_{jh}\)- \hat r_{2m}^j - \l
y_m^je^{2\l \phi_m^j} -\hat
r_m^j\]\mE\(\frac{p_m^{j+1} - 2p_m^j +
p_m^{j-1}}{h^2}\ \!\Big|\ \!\cF_{jh}\)dx.
\end{array}
\end{equation}
Using $\hat z_m^0=\hat z_m^m=p_m^0=p_m^m=0$
again, we get that
\begin{eqnarray*}
\begin{array}{ll}
\ds\mE\sum_{j=1}^{m-1}\int_G\mE\(\frac {\hat
z_m^{j+1}-2\hat z_m^j +\hat z_m^{j-1}}{
h^2}\ \!\Big|\ \!\cF_{jh}\)\mE\(\frac{p_m^{j+1} -
2p_m^j + p_m^{j-1}}{h^2}\ \!\Big|\ \!\cF_{jh}\)dx\\
\ns \ds=\mE\sum_{j=2}^{m-1}\!\int_G\!\!
\mE(\hat z_m^j\ \!|\ \!\cF_{(j-1)h})\mE\(\frac{p_m^j\!-\!2p_m^{j-1}\!+\!p_m^{j-2}}{
h^4}\ \!\Big|\ \!\cF_{(j-1)h}\)dx
\\
\ns \ds\q\!-\!2\mE\!\sum_{j=1}^{m-1}\!\int_G\!\!
\hat z_m^j\mE\(\frac{p_m^{j+1}\! -\! 2p_m^j
\!+\!
p_m^{j-1}}{h^2}\ \!\Big|\ \!\cF_{jh}\)dx+\mE\!\sum_{j=1}^{m-2}\!\int_G\!
\hat
z_m^j\mE\(\frac{p_m^{j+2}\!-\!2p_m^{j+1}\!+p_m^j}{
h^4}\ \!\Big|\ \!\cF_{(j+1)h}\)dx\\
\ns \ds=\mE\int_G \hat
z_m^1\frac{\mE(p_m^3\ \!|\ \!\cF_{h})
-4\mE(p_m^2\ \!|\ \!\cF_{h})+5p_m^1}{
h^4}dx\\
\ns\ds\q  +\mE \int_G \hat
z_m^{m-1}\frac{4p_m^{m-1} +
\mE(p_m^{m-1}\ \!|\ \!\cF_{(m-2)h}) - 4p_m^{m-2} +
p_m^{m-3}}{ h^4}dx
\\
\ns \ds\q +\mE \sum_{j=2}^{m-2} \int_G
\frac{\mE(p_m^{j+2}\ \!|\ \!\cF_{jh})
-4\mE(p_m^{j+1}\ \!|\ \!\cF_{jh})+5p_m^j+\mE(p_m^j\ \!|\ \!\cF_{(j-1)h})-4p_m^{j-1}+p_m^{j-2}}{
h^4}\hat
z_m^jdx\\
\ns \ds=\mE\sum_{j=1}^{m-1}\int_G\hat
z_m^j\Big\{\sum_{j_1,j_2=1}^n\pa_{x_{j_2}}\[b^{j_1j_2}\pa_{x_{j_1}}
\mE\(\frac{p_m^{j+1} -2p_m^j + p_m^{j-1} }{
h^2}\ \!\Big|\ \!\cF_{jh}\)\]\\
\ns \ds\q\q\q\q\q\q\q-\frac{\mE(\hat
z_m^{j+1}\ \!|\ \!\cF_{jh}) e^{-2\l \phi_m^{j+1}}-2\hat
z_m^j e^{-2\l \phi_m^j}+\hat z_m^{j-1} e^{-2\l
\phi_m^{j-1}}}{ h^2}\Big\}dx.
\end{array}
\end{eqnarray*}
Noting that
$z_m^j|_{\G}=p_m^j|_{\G}=0$ $(j=0, 1, \cdots, m)$,  one has
\begin{eqnarray*}
\begin{array}{ll} \ds
\mE\sum_{j=1}^{m-1}\int_G
\sum_{j_1,j_2=1}^n\pa_{x_{j_2}}\big(b^{j_1j_2}\pa_{x_{j_1}}\hat
z_m^j\big)
\mE\(\frac{p_m^{j+1}-2p_m^j+p_m^{j-1}}{
h^2}\ \!\Big|\ \!\cF_{jh}\)dx\\
\ns \ds= \mE\sum_{j=1}^{m-1}\int_G \hat z_m^j
\sum_{j_1,j_2=1}^n\pa_{x_{j_2}}\(b^{j_1j_2}\pa_{x_{j_1}}\mE\(\frac{p_m^{j+1}-2p_m^j+p_m^{j-1}}{
h^2}\ \!\Big|\ \!\cF_{jh}\)\)dx.
\end{array}
\end{eqnarray*}
Then by  \eqref{cs33}, we obtain
that
\begin{equation}\label{cscI2}
\begin{array}{ll}
\ds0 =-\mE \sum_{j=1}^{m-1}\int_G\[\hat
z_m^j\frac{\mE(\hat z_m^{j+1}\ \!|\ \!\cF_{jh}) e^{-2\l
\phi_m^{j+1}}-2\hat z_m^j e^{-2\l \phi_m^j}+\hat
z_m^{j-1} e^{-2\l
\phi_m^{j-1}}}{h^2}\\
\ns \ds\qq + \(\frac{\mE(\hat
r_{1m}^{j+1}\ \!|\ \!\cF_{jh})-\hat r_{1m}^j}{h}+\hat
r_{2m}^j +\l y_m^je^{2\l \phi_m^j}+\hat
r_m^j\)\mE\(\frac{p_m^{j+1}-2p_m^j+p_m^{j-1}}{h^2}\ \!\Big|\ \!\cF_{jh}\)\]dx.
\end{array}
\end{equation}

It follows from Lemma \ref{5.24-lm2} that
\begin{equation}\label{cscI1}
\begin{array}{ll}
\ds -\mE \sum_{j=1}^{m-1} \int_G \[\hat
z_m^j\frac{\mE(\hat z_m^{j+1}\ \!|\ \!\cF_{jh})e^{-2\l
\phi_m^{j+1}} - 2\hat z_m^je^{-2\l \phi_m^j} +
\hat z_m^{j-1}e^{-2\l
\phi_m^{j-1}}}{h^2}\]dx \\
\ns \ds=\mE\sum_{j=0}^{m-1}
\int_G\Big\{\frac{(\hat z_m^{j+1}-\hat
z_m^j)}{h}\frac{(\hat z_m^{j+1}e^{-2\l
\phi_m^{j+1}}-\hat
z_m^je^{-2\l \phi_m^j})}{ h} \\
\ns \ds=\mE\sum_{j=0}^{m-1} \int_G\[\frac{(\hat
z_m^{j+1}\!-\!\hat z_m^j)^2}{h^2}e^{-2\l
\phi_m^j}\!+\!\frac{\hat z_m^{j+1}\!-\!\hat
z_m^j}{h}\frac{e^{-2\l \phi_m^{j+1}}\!-\!e^{-2\l
\phi_m^j}}{h}\hat z_m^{j+1}\]dx,
\end{array}
\end{equation}
By Lemma \ref{5.24-lm2} and $p_m^j=K\hat r_m^j$,
we have that
\begin{equation}\label{cscI1.1}
\begin{array}{ll}
\ds -\mE \sum_{j=1}^{m-1} \int_G \[ \hat
r_m^j\mE\(\frac{p_m^{j+1} - 2p_m^j +
p_m^{j-1}}{h^2}\ \!\Big|\ \!\cF_{jh}\)\]dx
=K\mE\sum_{j=0}^{m-1} \int_G\frac{(\hat
r_m^{j+1}- \hat r_m^j)^2}{h^2}dx.
\end{array}
\end{equation}
Further, by (\ref{5.23-eq4}) and
Lemma \ref{5.24-lm2}, we get that
\begin{equation}\label{csp1}
\begin{array}{ll}
\ds
-\mE\sum_{j=1}^{m-1}\int_G\[\mE\(\frac{\hat
r_{1m}^{j+1}-\hat r_{1m}^j}{ h}\ \!\Big|\ \!\cF_{jh}\)
\]\mE\(\frac{p_m^{j+1}-2p_m^j+p_m^{j-1}}{
h^2}\ \!\Big|\ \!\cF_{jh}\)dx\\
\ns \ds=-\mE\sum_{j=1}^{m-1}\int_G
\mE\(\frac{\hat r_{1m}^{j+1}-\hat r_{1m}^j}{
h}\ \!\Big|\ \!\cF_{jh}\) \frac {1}{
h}\mE\(\frac{p_m^{j+1}-p_m^j}{ h}
-\frac{p_m^j-p_m^{j-1}}{
h}\ \!\Big|\ \!\cF_{jh}\)dx\\
\ns
\ds=\mE\sum_{j=1}^{m-1}\int_G\frac{\varrho}{\l^2}\mE\(\frac{\hat
r_{1m}^{j+1}-\hat r_{1m}^j}{
h}\ \!\Big|\ \!\cF_{jh}\)\mE\(\frac {\hat
r_{1m}^{j+1}e^{-2\l \phi_m^{j+1}}-\hat
r_{1m}^je^{-2\l \phi_m^j}}{
h}\ \!\Big|\ \!\cF_{jh}\)dx \\
\ns \ds= \mE \sum_{j=1}^{m-1} \int_G
\frac{\varrho}{\l^2}\[\frac{[\mE(\hat
r_{1m}^{j+1} - \hat r_{1m}^j\ \!|\ \!\cF_{jh})]^2}{
h^2}e^{-2\l \phi_m^j}\\
\ns\ds\qq + \mE\(\frac{\hat r_{1m}^{j+1} - \hat
r_{1m}^j}{h}\ \!\Big|\ \!\cF_{jh}\) \frac{e^{-2\l
\phi_m^{j+1}} - e^{-2\l \phi_m^j}}{ h}\mE(\hat
r_{1m}^{j+1}\ \!|\ \!\cF_{jh})\] dx,
\end{array}
\end{equation}
and
\begin{equation}\label{csp1.1}
\begin{array}{ll}
\ds -\mE\sum_{j=1}^{m-1}\int_G\[\l y_m^je^{2\l
\phi_m^j}+\hat r_{2m}^j
\]\mE\(\frac{p_m^{j+1}-2p_m^j+p_m^{j-1}}{
h^2}\ \!\Big|\ \!\cF_{jh}\)dx\\
\ns \ds=-\mE\sum_{j=1}^{m-1}\int_G \l
y_m^je^{2\l \phi_m^j} \frac {1}{
h}\mE\(\frac{p_m^{j+1}-p_m^j}{ h}
-\frac{p_m^j-p_m^{j-1}}{
h}\ \!\Big|\ \!\cF_{jh}\)dx\\
\ns
\ds\q+\mE\sum_{j=0}^{m-1}\int_G \frac{(\hat r_{2m}^{j+1}-\hat r_{2m}^j)}{h}\frac{p_m^{j+1}-p_m^j}{ h}dx\\
\ns
\ds=\mE\sum_{j=1}^{m-1}\int_G\frac{\varrho}{\l^2}\l
y_m^je^{2\l \phi_m^j}\mE\(\frac {\hat
r_{1m}^{j+1}e^{-2\l \phi_m^{j+1}}-\hat
r_{1m}^je^{-2\l \phi_m^j}}{
h}\ \!\Big|\ \!\cF_{jh}\)dx\\
\ns \ds\q+\mE\sum_{j=0}^{m-1}\int_G
\frac{\varrho}{\l^4}\frac{(\hat
r_{2m}^{j+1}-\hat r_{2m}^j)}{h} \frac{(\hat
r_{2m}^{j+1}e^{-2\l \phi_m^{j+1}}-\hat
r_{2m}^je^{-2\l \phi_m^j})}{
h}dx\\
\ns
\ds=\l\mE\sum_{j=1}^{m-1}\int_G\frac{\varrho}{\l^2}y_m^j\[
\mE\(\frac{\hat r_{1m}^{j+1}-\hat r_{1m}^j}{ h}
e^{-2\l \phi_m^j}\ \!\Big|\ \!\cF_{jh}\)\\
\ns\ds\quad+\frac{(e^{-2\l
\phi_m^{j+1}}\!\!-\!e^{-2\l \phi_m^j})}{
h}\mE(\hat r_{1m}^{j+1}\ \!|\ \!\cF_{jh})\]dx\\
\ns \ds\q+\mE\sum_{j=0}^{m-1}\int_G
\frac{\varrho}{\l^4}\[\frac{(\hat
r_{2m}^{j+1}-\hat r_{2m}^j)^2}{ h^2}e^{-2\l
\phi_m^j}  +\frac{(\hat r_{2m}^{j+1}-\hat
r_{2m}^j)}{ h} \frac{(e^{-2\l
\phi_m^{j+1}}-e^{-2\l \phi_m^j})}{ h}\hat
r_{2m}^{j+1}\]dx.
\end{array}
\end{equation}
By (\ref{cscI1})-(\ref{csp1.1}), it follows that
\begin{eqnarray*}
\begin{array}{ll}
\ds\mE\sum_{j=0}^{m-1} \int_G\[\frac{(\hat
z_m^{j+1}-\hat z_m^j)^2}{h^2}e^{-2\l
\phi_m^j}+\frac{\varrho}{\l^2}\frac{[\mE(\hat
r_{1m}^{j+1}-\hat
r_{1m}^j\ \!|\ \!\cF_{jh})]^2}{h^2}e^{-2\l \phi_m^j}\\
\ns\ds\qq\qq\q+\frac{\varrho}{\l^4}\frac{(\hat
r_{2m}^{j+1}-\hat r_{2m}^j)^2}{h^2}e^{-2\l \phi_m^j}+K\frac{\big(\hat r_m^{j+1}-\hat r_m^j\big)^2}{ h^2}\]dx\\
\ns\ds= -\mE\sum_{j=0}^{m-1} \int_G\frac{(\hat
z_m^{j+1}-\hat z_m^j)}{h}\frac{(e^{-2\l
\phi_m^{j+1}}-e^{-2\l \phi_m^j})}{h}\hat
z_m^{j+1}dx\\
\ns
\ds\q-\mE\sum_{j=1}^{m-1}\int_G\frac{\varrho}{\l^2}\mE\(\frac{\hat
r_{1m}^{j+1}-\hat r_{1m}^j}{h}\ \!\Big|\ \!\cF_{jh}\)
\frac{(e^{-2\l \phi_m^{j+1}}-e^{-2\l \phi_m^j})}{h}\mE(\hat r_{1m}^{j+1}\ \!|\ \!\cF_{jh})dx\\
\ns\ds\q-\l\mE\sum_{j=1}^{m-1}\int_G\frac{\varrho}{\l^2}y_m^j\[
\mE\(\frac{\hat r_{1m}^{j+1}-\hat r_{1m}^j}{ h}
e^{-2\l \phi_m^j}\ \!\Big|\ \!\cF_{jh}\)+\frac{(e^{-2\l
\phi_m^{j+1}}\!\!-\!e^{-2\l \phi_m^j})}{
h}\mE(\hat r_{1m}^{j+1}\ \!|\ \!\cF_{jh})\]dx\\
\ns\ds\q-\mE\sum_{j=0}^{m-1}\int_G
\frac{\varrho}{\l^4}\frac{(\hat
r_{2m}^{j+1}-\hat r_{2m}^j)}{h}\frac {(e^{-2\l
\phi_m^{j+1}}-e^{-2\l \phi_m^j})}{ h}\hat
r_{2m}^{j+1}dx.
\end{array}
\end{eqnarray*}
By  H\"older's inequality and the above equality,
 there is a positive constant
$C=C(K,\l)$, independent of $m$, such that
\begin{equation}\label{llK1}
\begin{array}{ll}
\ds\mE\sum_{j=0}^{m-1} \int_G\Big\{\frac{(\hat
z_m^{j+1}-\hat z_m^j)^2}{h^2}e^{-2\l
\phi_m^j}+\frac{\varrho}{\l^2}\frac{[\mE(\hat
r_{1m}^{j+1}-\hat
r_{1m}^j\ \!|\ \!\cF_{jh})]^2}{ h^2}e^{-2\l \phi_m^j}\\
\ns \ds\qq\qq+\frac{\varrho}{\l^4}\frac{(\hat
r_{2m}^{j+1}-\hat r_{2m}^j)^2}{
h^2}e^{-2\l \phi_m^j}+K\frac{\big(\hat r_m^{j+1} - \hat r_m^j\big)^2}{ h^2}\Big\}dx\\
\ns \ds\le C\mE\[\sum_{j=1}^{m-1}\int_G\(|\hat
z_m^j|^2+|\hat r_{1m}^j|^2+|\hat r_{2m}^j|^2
+K|\hat r_m^j|^2+|y_m^j|^2\)dx  +\int_G |\hat
r_{1m}^m|^2dx\].
\end{array}
\end{equation}
Finally,  by (\ref{llK1}) and
(\ref{5.23-eq6}),  recalling that $y\in
L^2_\mathbb{F}(\O;C([0,T];L^2(G)))$, we
get the desired estimate (\ref{5.24-eq7}).
This completes the proof of Proposition
\ref{5.24-prop1}.
\endpf

%%%%%%%%%%%%%%%%%%%%%%%%%%%%%%%%%%%%%%%%%%%%%%%%%%%%%

\section{Appendix B: Proof of \eqref{5.26-eq3}}

%%%%%%%%%%%%%%%%%%%%%%%%%%%%%%%%%%%%%%%%%%%%%%%%%%%%%%%%%

This appendix is addressed to proving \eqref{5.26-eq3}.

By
\eqref{csz7.1}, for a.e. $\o\in
\O$, $z_\o = \tilde z(\o)\in H^1(0, T; L^2(G))$  is a weak
solution to the following random equation:
\begin{equation}\label{5.26-eq20}
\left\{
\begin{array}{ll}
\ds \cA  z_\o=\tilde  r_{1,\o,t}+\tilde  r_{2,\o}+\l y_\o e^{2\l \phi}+\check  r_\o &\mbox{ in } Q, \\[3mm]
\ns\ds  z_\o=0   &\mbox{ on } \Si,\\[3mm]
\ns\ds   z_\o(0)=z_\o(T)=0  &\mbox{
in }G.
\end{array}
\right.
\end{equation}
Here $\tilde  r_{1,\o,t}= \tilde  r_{1,t}(\o)$,
$\tilde  r_{2,\o} = \tilde  r_{2}(\o)$, $y_{\o}
= y(\o)$ and $\check r_{\o} = \check r(\o)$.
Also, set $h_\o=\tilde  r_{1,\o,t}+\tilde  r_{2,\o}+\l y_\o e^{2\l \phi}+\check  r_\o$.

In the following, without loss of generality, we assume that $z_{\o}$ is smooth and give a uniform estimate for it.
Let $0<t_1<t_2<T$.  Multiplying the first
equation of \eqref{5.26-eq20} by
$t^2(T-t)^2z_{\o}$ and  integrating it in $(0,T)
\times G$, we
 get that
\begin{equation}\label{5.27-eq3}
\int_{t_1}^{t_2}\int_G  |\nabla z_{\o}|^2
 dxdt \leq C\int_{0}^{T}\int_G
\big(|h_{\o}|^2 + |z_{\o,t}|^2\big)dxdt.
\end{equation}
Put
$$
E(t) =\frac{1}{2}\int_G \big( |z_{\o,t}(t)|^2 +
|\nabla z_{\o}(t)|^2 \big)dx.
$$
By the usual energy estimate for the first
equation of \eqref{5.26-eq20} and noting the
time reversibility of \eqref{5.26-eq20}, we have
that
\begin{equation}\label{5.27-eq4}
E(t)\leq C\[E(s) +
\int_{t_1}^{t_2}\int_G|h_{\o}(\tau,x)|^2dxd\tau\],
\q\forall\  t, s \in [t_1,t_2].
\end{equation}

Integrating \eqref{5.27-eq4} with respect to $s$
from $t_1$ to $t_2$, we obtain
\begin{equation}\label{5.27-eq5}
E(t)\leq C\[\int_{t_1}^{t_2} E(s)ds +
\int_{t_1}^{t_2}\int_G|h_{\o}(\tau,x)|^2dxd\tau\],
\q\forall\ t \in [t_1,t_2].
\end{equation}
By \eqref{5.27-eq3} and \eqref{5.27-eq5},  for any  $t \in [t_1,t_2]$, there
is a constant $C>0$  such that
\begin{equation}\label{5.27-eq10}
|z_{\o,t}(t)|^2_{L^2(G)} +
|z_{\o}(t)|^2_{H^1_0(G)}\leq C\int_{0}^{T}\int_G
\big(|h_{\o}|^2 + |z_{\o,t}|^2\big)dxdt.
\end{equation}

Applying the usual energy estimate to the first
equation of \eqref{5.26-eq20} and  noting
the time reversibility of \eqref{5.26-eq20} again,
similar to the proof of \eqref{5.27-eq4}, we find that
$$
\begin{array}{ll}\ds
|z_{\o} |^2_{C([0,T];H^1_0(G))\cap
C^1([0,T];L^2(G))}\\
\ns\ds \leq C\[  |z_{\o,t}(t)|^2_{L^2(G)} +
|z_{\o}(t)|^2_{H^1_0(G)}
+|h_{\o}|^2_{L^2(0,T;L^2(G))}\].
\end{array}
$$
This, together with \eqref{5.27-eq10}, implies
that
$$
|z_{\o} |^2_{C([0,T];H^1_0(G))\cap
C^1([0,T];L^2(G))} \leq C\int_{0}^{T}\int_G
\big(|h_{\o}|^2 + |z_{\o,t}|^2\big)dxdt.
$$
It follows  that
$$
\mE|z_\o|_{C([0,T];H_0^1(\O)))\cap
C^1([0,T];L^2(\O))}^2 \leq
C\mE\int_{0}^{T}\int_G \big(|h_{\o}|^2 +
|z_{\o,t}|^2\big)dxdt.
$$
This, together with $\tilde z \in
L^2_\dbF(\O;H^1(0,T;L^2(G)))$, implies that
$$
\tilde z \in
L^2_\dbF(\Omega;C([0,T];H_0^1(\O)))\cap
L^2_\dbF(\Omega;C^1([0,T];L^2(\O))).
$$
\endpf

%%%%%%%%%%%%%%%%%%%%%%%%%%%%%%%%%%%%%%%%%%%%%%%%%%%%%%%5

\end{document}